\documentclass{amsart}
\usepackage[utf8]{inputenc}

\usepackage{amsmath}
\usepackage{amsthm}
\usepackage{amsfonts}
\usepackage{amssymb}
\usepackage{mathabx}
\usepackage{mathtools}
\usepackage{mdframed}
\usepackage{hyperref}
\usepackage{multicol}
\usepackage{tikz-cd}
\usepackage{svg}
\usepackage{listings}

\newtheorem{theorem}{Theorem}
\newtheorem{lemma}[theorem]{Lemma}
\newtheorem{corollary}[theorem]{Corollary}
\newtheorem{proposition}[theorem]{Proposition}
\newtheorem{claim}{Claim}
\newtheorem{question}{Question}

\theoremstyle{definition}
\newtheorem*{remark}{Remark}
\newtheorem{definition}{Definition}

\usepackage{multirow}
\usepackage{multicol}
\usepackage{subfig}
\usepackage{graphicx}
\usepackage{wrapfig,lipsum}
\usepackage[top=1in,bottom=1in,left=1in,right=1in]{geometry} 

\newcommand{\orb}{(O,\psi)}
\newcommand{\orbprm}{(O,\psi')}
\newcommand{\orbone}{(O_1,\psi_1)}
\newcommand{\orbtwo}{(O_2,\psi_2)}
\newcommand{\orbnot}{(O_0,\psi_0)}
\newcommand{\orbi}{(O_i,\psi_i)}

\newcommand{\order}{\operatorname{order}}
\newcommand{\sym}{\operatorname{Sym}}
\newcommand{\hyp}{\operatorname{hyp}_c}
\newcommand{\MCG}{\operatorname{MCG}}
\newcommand{\cone}{\operatorname{cone}}
\newcommand{\HP}{\mathbb{H}}
\newcommand{\R}{\mathbb{R}}
\newcommand{\PSL}{\operatorname{PLS}}
\newcommand{\pld}{\operatorname{PLD}}
\newcommand{\G}{\mathcal{G}}
\newcommand{\bt}{\mathbf{t}}

\title{Flexible hyperbolic  cone metrics on the genus 2 surface}

\author{Katherine Chui}
\address{Department of Mathematics, Rice University, Houston, TX, USA}
\email{kac20@rice.edu}

\author{Jacob Russell}
\address{Department of Mathematics and Statistics, Swarthmore College, Swarthmore, PA, USA}
\email{jrussel2@swarthmore.edu}

\begin{document}
	
	\begin{abstract}
	A negatively curved hyperbolic cone metric on a surface is \emph{rigid} if it is determined by the support of its Liouville current. We use a theorem of Erlandsson, Leininger, and Sadanand to show that there are nine mapping class group orbits of equivalence classes of non-rigid (aka flexible) metrics in the case of the genus 2 surface.
\end{abstract}

\maketitle

\section{Introduction}
A cornerstone problem in the geometry of surfaces is identifying what data will determine the metric on the surface. The  \emph{Liouville geodesic current} has emerged as a central object in these questions, having been shown to determine the metric is a  variety of settings \cite{otal_length_rigidity,croke_npc_length_rigidity, HP_marked_length, frazier_length_rigidity,constantine_marked_length_rigidity}. In the cases of flat metrics and negatively curved hyperbolic  cone metrics, one does not usually need the full Liouville current  as the support of the current determines the metric in the vast majority of cases \cite{DELS_flat_metrics,ELS_hyp_cone_metrics}.

In the current paper, we focus on the case of negatively curved hyperbolic cone metrics on a closed surface $S$ of genus 2. We denote the isotopy classes of such metrics by $\hyp(S)$ and for each $\varphi \in \hyp(S)$, we let $\G_\varphi$ denote the support of the Liouville current for $\varphi$. We say $\varphi$ is \emph{rigid} if $\G_\varphi = \G_{\varphi'}$ implies $\varphi = \varphi'$. Otherwise, we say $\varphi$ is \emph{flexible}. If $\varphi$ is a flexible metric, the set of metrics $\varphi' \in \hyp(S)$ with $\G_{\varphi'} =\G_\varphi$ is called the \emph{flexibility class} of $\varphi$.

We prove  there are nine mapping class group orbits of flexibility classes in $\hyp(S)$ when $S$ has genus 2.

\begin{theorem}\label{intro_thm:42}
	If $S$ is the closed, connected, orientable surface of genus 2, then $\hyp(S)$ contains nine $\MCG(S)$-orbits of flexibility classes.
\end{theorem}

Our proof of Theorem \ref{intro_thm:42} rests on Erlandsson, Leininger, and Sadanand's characterization of flexible metrics in $\hyp(S)$ \cite{ELS_hyp_cone_metrics}. They prove that $\varphi \in \hyp(S)$ is flexible if and only if  $(S,\varphi)$ admits a locally isometric, finite degree branched cover of a non-triangular hyperbolic orbifold $\orb$ where every cone point is sent to an even order orbifold point. We call such a covering a \emph{flexible covering} by $(S,\varphi)$.  Our proof of Theorem \ref{intro_thm:42} is thus accomplished by enumerating the flexible coverings by the genus 2 surface.

The first step in our proof is to  create some  obstructions to having a flexible covering, and then produce the small list of degrees and orbifolds that avoid these obstructions; Sections \ref{section:basic_lemmas} and \ref{Section: list of Orbifolds}. Once we have this list, we determine what restrictions the existence of a flexible covering would impose on the cone points in $(S,\varphi)$ and on the local degrees of the preimage of orbifold points; Section \ref{sec:local_degree_restrictions}. 
This creates a very limited list of possible orbifolds that are covered and limits the cone points of any flexible metric.

\begin{theorem}\label{intro_thm:orbifolds_and_degrees}
	Let $S$ be a closed, orientable surface of genus $2$. If $p\colon (S,\varphi) \to \orb$ is a finite degree flexible covering, then we have one of the following for the degree of $p$, the signature of $\orb$, and the cone points of $(S,\varphi)$:
	\begin{center}
		\begin{tabular}{c|c|c}
			Degree of $p$ & Signature of $\orb$ & Cone points in $(S,\varphi)$  \\ \hline
			$3$ & $(0;2,2,2,3)$ & 1 point with angle $3\pi$ \\ 
			$3$ & $(0;2,2,3,3)$ & 2 points with angle $3\pi$ \\ 
			$3$ & $(0;2,2,2,3)$ & 3 points with angle $3\pi$ \\ 
			$3$ & $(1;2)$ & 1 point with angle $3\pi$  \\ 
			$4$ & $(0;2,2,2,4)$ & 1 point with angle $4\pi$  \\
			$6$ & $(0;2,2,2,4)$ & 1 point with angle $3\pi$  \\
			$6$ & $(0;2,2,2,3)$ & 2 points with angle $3\pi$ or 1 point with angle $4\pi$  \\
		\end{tabular}
	\end{center}
	In particular, if $(S,\varphi)$ has three or more cone points, then $\varphi$ is rigid.
\end{theorem}

The second step of our proof is showing that in each of the cases listed in Theorem \ref{intro_thm:orbifolds_and_degrees}, there does exist a flexible covering map.  This is accomplished using the correspondence between finite covers and homomorphism of  $\pi_1(O - \{\text{orbifold points}\})$ into the symmetric group; see Section \ref{subsec:permutations}. This correspondence allows us to reduce the problem of finding flexible covers to finding tuples of elements of the symmetric group with specific cycle structures and satisfying some simple equations. Since the degrees of our covers are at most $6$, we can use computer code to perform an exhaustive search of the all tuples in the symmetric group satisfying these conditions. This creates a list of all the covering space equivalence classes of flexible covering maps.

The final step of our proof is to count the distinct mapping class groups orbits of flexibility classes coming from this finite list of covering maps. The challenge here is that topologically distinct covering maps can give rise to the same orbit of flexibility classes in $\hyp(S)$.  This over counting is resolved by using the work of Erlandsson, Leininger, and Sadanand to show that two flexible metric are in the same mapping class group orbit if and only if they admit flexible coverings of a canonical orbifold that differ by  a homeomorphism of that orbifold; see Section \ref{sec:flexible_cone_metrics}. 

Erlandsson, Leininger, and Sadanand originally showed that for a surface of genus $g$, any metric with at least $32(g-1)$ cone points is rigid. Our calculations show that in the  $g=2$ case, just 3 cone points are sufficient to ensure rigidity.  Erlandsson, Leininger, and Sadanand note that you can always find a flexible metric with $g-1$ cone points, so this begs the following question.

\begin{question}
	What is the smallest positive integer $k$ so that for all genus $g$,  if $\varphi \in \hyp(S)$ has $k(g-1)$ cone points, then $\varphi$ is rigid?
\end{question}

\noindent\textbf{Acknowledgments:} The authors are grateful to Chris Leininger for both suggesting the topic and for providing guidance through out the course of the project. Russell was supported by NSF grant DMS-2103191.

\section{Preliminaries}
Throughout, we will only consider connected, closed, and orientable  surfaces with finite genus. A \emph{hyperbolic cone metric} on a surface  is a metric that is locally isometric to the hyperbolic plane except for a finite number of  cone singularities. We call the cone singularities the \emph{cone points} of the metric and the remaining points the \emph{regular points}. We will assume that our hyperbolic cone metric always contain at least one cone point. We use the notation $(S, \varphi)$ to denote the topological surface $S$ equipped with the hyperbolic cone metric $\varphi$.

Given a hyperbolic cone metric $\varphi$ on $S$, we let $\cone(\varphi)$ denote the list of cone angles for the cone points of $\varphi$ with multiplicity. For each  point $x$  of $(S,\varphi)$ we use $\Theta_x$ to denote the  angle of $x$. We say a hyperbolic cone metric is \emph{negatively curved} if the cone angle of each cone point is strictly larger than $2\pi$. We let $\hyp(S)$  denote the set of isotopy classes of negatively curved hyperbolic cone metrics on the surface $S$. By abuse of notation, we will not distinguish between a metric and its isotopy class.

Homeomorphisms of $S$ act on $\hyp(S)$ by pushing forward the metrics. Since the metrics in $\hyp(S)$ are considered up to isotopy, a homeomorphism that is isotopic to the identity will act by the identity on $\hyp(S)$. Thus, we get an action of the \emph{mapping class group} $\MCG(S)$ on $\hyp(S)$, where $\MCG(S)$ is the group of homeomorphism of $S$ up to isotopy.

A \emph{hyperbolic orbifold} $\orb$ is a  surface $O$ equipped with a hyperbolic cone metric $\psi$ where each the cone angles is of the form $2\pi / r$ for some integer $r \ge 2$.   We often call the cone points of an orbifold the \emph{orbifold points}.   If an orbifold point $x$ has cone angle $\Theta_x = 2\pi /r$, then we say $x$ has \emph{order} $r$, and use  $\order(x)$ to denote this order. The \emph{signature} of an orbifold is a tuple $(g;r_1,\dots,r_m)$ where $g$ is the genus of the underlining topological surface, $m$ is the number of orbifold points, and the $r_i$ are the orders of the orbifold points. As a convention, we will always arrange the orders in the signature from lowest to highest. An \emph{orbifold homeomorphism}, $f \colon \orb \to \orbprm$ is a homeomorphism of $O$ that sends each orbifold point to an orbifold point of the same order. The \emph{orbifold mapping class group}, $\MCG\orb$, is the group of these orbifold homeomorphisms up to isotopy. By abuse of notation, when we say $f\colon \orb \to \orbprm$ is a homeomorphism, we mean that it is an orbifold homeomorphism. 

\subsection{Area, covering maps, and local degrees} 
\label{Gauss-Bonnet}

Let $S$ be a surface of genus $g \geq 0$ equipped with a hyperbolic cone metric $\varphi$ ($\varphi$ is not necessarily negatively curved). If $\cone(\varphi) = (\Theta_1,\dots,\Theta_n)$, then the \emph{area} of $\varphi$, $A(\varphi)$, is given by the Gauss--Bonnet theorem: $$A(\varphi) =  -2\pi\left(2-2g-\sum_{i=1}^n\left(1-\frac{\Theta_i}{2\pi}\right)\right).$$
In the case of when $(S,\varphi)$ is an orbifold $\orb$ with signature $(g;r_1,\dots,r_m)$, this gives $$A(\psi) = -2\pi\left(2-2g-\sum_{i=1}^m\left(1-\frac{1}{r_i}\right)\right).$$

Fix $\varphi \in \hyp(S)$ and a hyperbolic orbifold $\orb$. Suppose there is a locally isometric, branched covering map  $p\colon (S,\varphi) \to \orb$ with finite degree $D \geq 1$. The following results will apply to any finite degree, locally isometric branched cover of hyperbolic cone metrics, but we restrict to these particular kinds of metrics as this will be the only situation where we will apply them.

The areas of $\varphi$ and $\psi$ are related by the \emph{Area Covering Formula}:
\begin{equation}\label{Area Covering Theorem}
	D = \frac{A(\varphi)}{A(\psi)}.
\end{equation}

For any $y \in (S,\varphi)$, we let $\deg(y)$ denote the local degree of $p$ at $y$.  For each orbifold point $x \in \orb$ and each $y \in p^{-1}(x)$ we have three basic facts from covering space theory:
\begin{equation}\label{eq:cover_degree}
	\Theta_x = \frac{\Theta_y}{\deg(y)},
\end{equation}
\begin{equation} \label{eq:sum_local_degree}
	D = \sum_{y\in p^{-1}(x)}\deg(y),
\end{equation}
\begin{equation} \label{eq:preimage_size}
	|p^{-1}(x)| = \frac{D}{\deg(y)}.  
\end{equation}

The local degrees of points in the preimage of a branched cover are going to play in an important role in the sequel. Thus, for any $x \in \orb$, we define $\pld(x)$ be the list (with multiplicity) of the local degrees of the points in $p^{-1}(x)$\footnote{PLD stands for ``preimage local degrees''.}.

\subsection{Puncturing surfaces,  developing maps, and holonomy}\label{subsec:punctures}
Let $\varphi \in \hyp(S)$, then define $\dot{S} = S - \cone(\varphi)$ to be the punctured surface obtained by deleting the cone points in $\varphi$. The restriction of $\varphi$ to $\dot{S}$ produces an incomplete hyperbolic metric $\dot{\varphi}$ on  $\dot{S}$. If $\widetilde{\dot{S}}$ is the universal cover of $\dot{S}$, let $(\hat{S}, 
\hat{\varphi})$ be the metric completion of $\widetilde{\dot{S}}$ with respect to the pull-back of the metric $\dot{\varphi}$. There is a $\pi_1(\dot{S})$-equivariant, locally isometric covering map $\hat{j} \colon \hat{S} \to (S,\varphi)$.

Since $\widetilde{\dot{S}}$ is simply connected and is equipped with a metric that is locally isometric to $\HP$, there is a \emph{developing map} $\widetilde{\dot{S}} \to \HP$. Up to post composing by an element of $\PSL(2,\R)$, this developing map is the unique orientation preserving local isometry of $\widetilde{\dot{S}}$ to $\HP$. The developing map extends to the metric completion  $\hat{S}$. We denote this developing map by $$ D\colon \hat{S} \to \HP.$$ 
As with the developing map  $\widetilde{\dot{S}} \to \HP$, if $D' \colon \hat{S} \to \HP$ is an orientation preserving local isometry, then there is $f \in \PSL(2,\R)$ so that $D = f \circ D '$.

Since $\pi_1(\dot{S})$ acts  by covering transformations on $\hat{S}$, for each $\gamma \in \pi_1(\dot{S})$, the map $D \circ \gamma \colon \hat{S} \to \HP$ is also an  orientation preserving local isometry. Thus, there is $\rho(\gamma) \in \PSL(2,\R)$ so that $D \circ \gamma = \rho(\gamma) \circ D$. Hence the map $\gamma \to \rho(\gamma)$ defines a homomorphism $\rho \colon \pi_1(\dot{S}) \to \PSL(2,R)$. This homomorphism is the \emph{holonomy representation} of $\pi_1(\dot{S})$ for the developing map $D$.

\subsection{Coverings via permutations}\label{subsec:permutations} 
A standard way of constructing  degree $D$ (non-branched) covering spaces is to use representations of the fundamental group into the finite  symmetric group $\sym(D)$; see, for example, \cite[Chapter 1]{Hatcer_book}. This construction characterizes covering spaces up to equivalence of covers.

\begin{theorem}\label{theorem:top_cover_permutaion}
	Let $\Sigma$ and $S$ be connected, compact,orientable surfaces, possibly with boundary. There is a one-to-one correspondence between
	\begin{enumerate}
		\item equivalence classes of degree $D$ covering map $p \colon S \to \Sigma$;
		\item conjugacy classes of homomorphism $\Phi \colon \pi_1(\Sigma) \to \sym(D)$ with transitive image.
	\end{enumerate}
\end{theorem}

\begin{remark}
	The requirement that the image of $\Phi$ is transitive is only to ensure that the covering space $S$ of $\Sigma$ is connected. Homomorphisms with non-transitive  images will correspond with disconnected covers. The same is true in Theorem \ref{theorem:banched_cover_permutation} below.
\end{remark}

There are several equivalent ways of constructing the permutation representation $\Phi_p \colon \pi_1(\Sigma) \to \sym(D)$ corresponding to the covering map $p\colon S \to \Sigma$. Here are two options that will be useful for future arguments.

\begin{enumerate}
	\item Fix a basepoint $x \in \Sigma$ for $\pi_1(\Sigma)$ and let $\{y_0,\dots y_{D-1}\} = p^{-1}(x)$. For each $\gamma \in \pi_1(\Sigma)$, $\gamma$ lifts to $D$ oriented paths in $S$ connecting the $y_i$.  The orientations on these paths plus the indices of the $y_i$ defines a permutation of the set $\{0,\dots,D-1\}$ (e.g. if $\gamma$ lifts to a path from $y_2$ to $y_4$ then this permutation sends 2 to 4). This permutation is $\Phi_p(\gamma)$. 
	\item Use the covering map $p$ to view $\pi_1(S)$ as a index $D$ subgroup of $\pi_1(\Sigma)$. $\Phi_p$ is then the homomorphism into $\sym(D)$ that comes from the standard action of $\pi_1(\Sigma)$ on the cosets of $\pi_1(S)$.
\end{enumerate}

In Theorem \ref{theorem:banched_cover_permutation} below, we present a modification of Theorem \ref{theorem:top_cover_permutaion} that allows one to construct branched covers with specific local degrees of the branch points.  For this we need a basic fact about extending degree $D$ covers of the circle to branched covers of closed disks.

\begin{lemma}\label{lemma:extending_boundary_maps}
	Fix a positive integer $D$. Let $Y$ and $X$ be closed disks. For every interior point $x \in X$ and every degree $D$ covering map of the circle $p \colon \partial Y \to \partial X$, there is  has an extension to a branched cover $p \colon Y \to X$ where the only branched point in $X$ is $x$ and  $p^{-1}(x)$ is a single point $y$ with $\deg(y) = D$.
\end{lemma}

Now, given a closed surface $\Sigma$ and a collection of points $x_1,\dots,x_m \in \Sigma$, let $\overline{\Sigma}$ denote the surface with boundary obtained by deleting a small open neighborhood of each point $x_i$ so that $\overline{\Sigma}$ is a surface with $m$ boundary component. We characterize branched covers of $\Sigma$ by representations of $\pi_1(\overline{\Sigma})$ into the symmetric group where the cycle structure of the images of the boundary curves correspond to the local degrees of the preimages of the branch points. The  idea is to apply Theorem \ref{theorem:top_cover_permutaion} to the representation of $\pi_1(\overline{\Sigma})$ to get a non-branched covering of $\overline{\Sigma}$ and then use Lemma \ref{lemma:extending_boundary_maps} to extend to a branched covering of $\Sigma$ with the correct local degrees.

\begin{theorem}\label{theorem:banched_cover_permutation}
	Let $S$ and $\Sigma$ be closed, connected  surfaces, and let $x_1,\dots, x_m$ be a set of points on $\Sigma$. Let $\delta_i$ be the boundary component of $\overline{\Sigma}$ corresponding to the  point $x_i$. Let $D$ be a positive integer and for each $i \in \{1,\dots,m\}$, let $\ell_1^i,\dots, \ell_{n_i}^i$ be a collection of positive integers so that $\sum_{j=1}^{n_i} \ell_j^i = D$.
	
	There is a one-to-one  correspondence between:
	\begin{enumerate}
		\item Equivalence classes of degree $D$ branched covers $p \colon S \to \Sigma$  where the only branched points are $x_1,\dots, x_m$ and $\pld(x_i) = (\ell_1^i,\dots,\ell_{n_i}^i)$ for each $i \in \{ 1,\dots,m\}$. 
		\item Conjugacy classes of homomorphisms $\Phi \colon \pi_1(\overline{\Sigma}) \to \sym(D)$ with transitive image  where for each $i \in \{1,\dots,m\}$ and for any element  $\gamma \in \pi_1(\overline{\Sigma})$ in the free homotopy class of $\delta_i$, $\Phi(\gamma)$ is a permutations with cycle structure $(\ell_1^i,\dots, \ell_{n_i}^i)$.
	\end{enumerate}
\end{theorem}

\begin{remark}
	If we  equip $\Sigma$ with a hyperbolic cone metric  so that it is a hyperbolic orbifold $\orb$ where $x_1,\dots, x_m$ are the orbifold points, then each homomorphism $\pi_1(\overline{O}) \to \sym(D)$ corresponds to a branched cover $p \colon S \to O$. We can pull back the metric on $\orb$ under $p$ to  equip $S$ with a hyperbolic metric $\varphi$  so that $p$ is a locally isometric branched cover $p \colon (S,\varphi) \to \orb$. Moreover, the angles of any cone points of $\varphi$ are determined by the cycle structure in $\sym(D)$ and the angles of the orbifold points on $\orb$. 
\end{remark}

\begin{proof}
	Assume there is degree $D$ branched cover $p \colon S \to \Sigma$ so that  the only branched points are $x_1,\dots,x_m$ and  $\pld(x_i) = (\ell_1^i,\dots,\ell_{n_i}^i)$ for each $i \in \{1,\dots,m\}$. Let $y_1^i,\dots, y_{n_i}^i$ be the preimages of $x_i$ under $p$ so that the local degree of $p$ at $y_j^i$ is $\ell_j^i$. Now, $\overline{S} = p^{-1}(\overline{\Sigma})$ is a subsurface of $S$ so that the restriction of $p$ to $\overline{S}$ is a  (non-branched) degree $D$ cover $\overline{p} \colon \overline{S} \to \overline{\Sigma}$ of surfaces with boundary. By Theorem \ref{theorem:top_cover_permutaion}, this cover corresponds to a homomorphism $\Phi \colon \pi_1(\overline{\Sigma}) \to \sym(D)$ with transitive image.
	
	For each $i \in \{1,\dots,m\}$, the preimage $\overline{p}^{-1}(\delta_i)$ is a set of boundary components $\eta_1^i,\dots, \eta_{n_i}^i$ of $\overline{S}$ where each $\eta_j^i$ in the boundary of a disk in $S$ containing $y_j^i$. The restriction of $\overline{p}$ to $\eta_j^i$ is a covering map of the circle with degree $\ell_j^i$, the local degree of $p$ at $y_j^i$. Hence, for any element in $\gamma \in \pi_1(\overline{\Sigma})$ in the free homotopy class of $\delta_i$, $\Phi(\gamma)$ is a permutation with cycle structure $(\ell_1^i,\dots, \ell_{n_i}^i)$.
	
	Now assume that there exists the desired homomorphism $\Phi \colon \pi_1(\overline{\Sigma}) \to \sym(D)$. By Theorem \ref{theorem:top_cover_permutaion}, there exists a corresponding degree $D$ cover $\overline{p}\colon \overline{S} \to \overline{\Sigma}$ where $\overline{S}$ is a connected surface with boundary.   For each $i \in \{1,\dots,m\}$, the preimage $\overline{p}^{-1}(\delta_i)$ is a set $\{\eta_1^i,\dots, \eta_{n_i}^i \}$ of boundary components of $\overline{S}$. Moreover, the restriction of $\overline{p}$ to $\eta_j^i$ is a degree $\ell_j^i$ cover of the circle.  Let $r_j^i$ denote this restriction of $\overline{p}$ to $\eta_j^i$.
	
	For each $i\in \{1,\dots,m\}$ and $j \in \{1,\dots,n_i\}$, let  $Y_j^i$ be a closed disk. Let $X_i$ be the closed disk on $\Sigma$ centered at $x_i$ whose interior is removed to create $\overline{\Sigma}$. Be identifying $\eta_j^i$ with $\partial Y_j^i$, we have a degree  $\ell_j^i$ covering map $r_j^i \colon \partial Y_j^i \to \partial x_i$.  By Lemma \ref{lemma:extending_boundary_maps}, this map can be extended to a degree $\ell_j^i$ branched cover $r_j^i \colon Y_j^i \to X_i$ where  ${r_j^i}^{-1}(x_i)$ is a single point $y_j^i$ with $\deg(y_j^i)=\ell_j^i$.
	
	We build the  surface $S$ and the branched cover $p \colon S \to \Sigma$ in two steps:
	\begin{enumerate}
		\item Create the closed surface $S$ by  attaching the boundary of the disk $Y_j^i$ to $\eta_j^i$ for each $i\in \{1,\dots,m\}$ and $j \in \{1,\dots,n_i\}$.
		\item Define $p\colon S \to \Sigma$ by $p(z) = \overline{p}(z)$ if $z \in \overline{S}$ and $p(z) = r_j^i(z)$ if $z \in Y_j^i$.
	\end{enumerate}
\end{proof}

\subsubsection{Permutation representations via tuples} \label{subsec:tuples}
We will apply Theorem \ref{theorem:banched_cover_permutation} when $\Sigma$ is either a sphere with four distinguished point $x_1,x_2,x_3,x_4$ or a torus with a single disguised point $x$. 

\begin{figure}[ht]
	\centering
	\def\svgwidth{0.35\columnwidth}
    \hspace{2cm}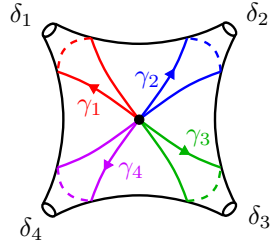
	\caption{Generators of $\pi_1(\overline{\Sigma})$ when $\Sigma$ is a sphere with four distinguished points.}%
	\label{fig:4-holed_sphere}%
\end{figure}

In the sphere case,  if $\gamma_1,\gamma_2,\gamma_3, \gamma_4$ are the loops shown in Figure \ref{fig:4-holed_sphere}, then  \[ \pi_1(\overline{\Sigma}) = \langle \gamma_1,\gamma_2,\gamma_3, \gamma_4 \mid \gamma_1\gamma_2\gamma_3 =\gamma_4^{-1} \rangle\] is a presentation for $ \pi_1(\overline{\Sigma})$. We can therefore specify a homomorphism $\Phi \colon \pi_1(\overline{\Sigma}) \to \sym(D)$ by defining each $\Phi(\gamma_i)$ to be $s_i \in \sym(D)$ so that $s_1s_2s_3= s_4^{-1}$. Thus, homomorphisms of $\pi_1(\overline{\Sigma)}$ to $\sym(D)$ are in one-to-on correspondence with tuples, $[s_1,s_2,s_3,s_4]$, of elements of $\sym(D)$ where $s_1s_2s_3= s_4^{-1}$. We say two such tuples, $[s_1,s_2,s_3,s_4]$ and $[s'_1, s'_2, s'_3, s'_4]$, are \emph{conjugate} if there exists a single $t \in \sym(D)$ so that for all $i \in  \{1,2,3,4\}$,  $ts_it^{-1} = s'_i$. Thus, two homomorphisms are conjugate if and only if their corresponding tuple are conjugate.

\begin{figure}[h]
	\begin{center}
		
        \hspace{2cm} \def\svgwidth{0.75\columnwidth}
        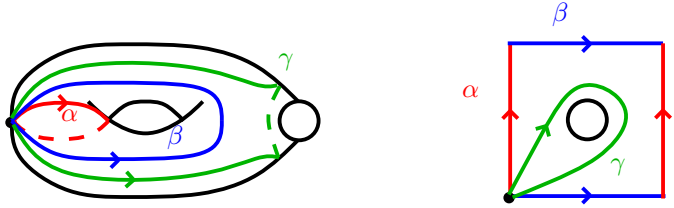
	\end{center}

	\caption{Generators for $\pi_1(\overline{\Sigma})$ when $\Sigma$ is a torus with a single distinguished point.}
	\label{fig:1-holed_torus}
\end{figure}

When $\Sigma$ is a torus with a single distinguished point $x$, then    $$ \pi_1(\overline{\Sigma}) = \langle \alpha,\beta, \gamma \mid \alpha \beta\alpha^{-1} \beta^{-1} = \gamma \rangle$$ where  $\alpha,\beta,\gamma$ are the loops shown in Figure \ref{fig:1-holed_torus}. As in the proceeding paragraph, there is a one-to-one correspondence between  tuples, $[s_1,s_2,s_3]$, of elements in $\sym(D)$ with $s_1s_2s_1^{-1}s_2^{-1} = s_3$ and homomorphisms $\Phi \colon \pi_1(\overline{\Sigma}) \to \sym(D)$  where  $$\Phi(\alpha) = s_1,\ \Phi(\beta)=s_2,\ \Phi(\gamma) = s_3.$$ We define  two tuples to be conjugate just as in the proceeding paragraph so that two tuples are conjugate if and only if the homomorphisms are conjugate.

\subsubsection{Permutation representations of nested covers}\label{subsectnested_covers}
Let $S$, $\Sigma$, $F$ be closed connected surfaces. Suppose you have the following commuting diagram of finite degree branched covering maps.
\begin{center}
	\begin{tikzcd}[ sep=small]
		S \arrow[dr, "r"]\arrow[to=3-1,"p"]  &  \\
		& F \arrow[ld, "q"]  \\
		\Sigma & 
	\end{tikzcd}
\end{center}

Let $\overline{\Sigma}$ be the surface with boundary obtained by deleting a small neighborhood of the branched points of $p$. Let $\overline{F} = q^{-1}(\overline{\Sigma})$ and $\overline{S} = p^{-1}(\overline{\Sigma}) = r^{-1}(\overline{F})$. Similarly, let $F_0$ be the surface obtained by deleting a small neighborhood of the branched points of $r$ and let $S_0 = r^{-1}(F_0)$. Note $\overline{F} \subseteq F_0$ and $\overline{S} \subseteq S_0$, but these inclusions might not be equalities. The maps $p$, $q$ and $r$ induce the following non-branched covers.

\begin{center}
	\begin{tikzcd}[ sep=small]
		\overline{S} \arrow[dr, "\overline{r}"]\arrow[to=3-1,"\overline{p}"] \arrow[rr,hook]&  & S_0 \arrow[d, "r_0"]\\
		& \overline{F} \arrow[ld, "\overline{q}"] \arrow[r,hook] & F_0 \\
		\overline{\Sigma} & 
	\end{tikzcd}
\end{center}

Since $\overline{p}$, $\overline{q}$, and $\overline{r}$ induce injections on the fundamental group we consider $\pi_1(\overline{\Sigma})$, $\pi_1(\overline{F})$ and $\pi_1(\overline{S})$ as the chain of subgroup $$\pi_1(\overline{S}) < \pi_1(\overline{F}) < \pi_1(\overline{\Sigma}). $$ Similarly, we use $r_0$ to consider $\pi_1(S_0)$ as subgroup of $\pi_1(F_0)$. For each of these non-branched covers let $\Phi_\ast$ denote the corresponding permutation representation, where $\ast \in \{\overline{p}, \overline{r}, r_0 \}$.

\begin{lemma}\label{lem:nested_covers}
	The subgroup $\Phi_{\overline{p}}(\pi_1(\overline{F}))$ in $\sym(\deg(p))$ is isomorphic to $\Phi_{r_0}(\pi_1(F_0))$ in $\sym(\deg(r))$.
\end{lemma}

\begin{proof}
	First recall that $\ker(\Phi_{\overline{p}})$ and $\ker(\Phi_{\overline{r}})$ are both equal to $\pi_1(\overline{S})$ and $\ker(r_0) = \pi_1(S_0)$. This is because the kernel of a permutation representation of a finite cover is exactly the loops in the base that lift to a collections of loops in the cover. This implies that $\Phi_{\overline{p}}(\pi_1(\overline{F}))$ is isomorphic to $\Phi_{\overline{r}}(\pi_1(\overline{F}))$. Thus, it suffices to show that $\Phi_{\overline{r}}(\pi_1(\overline{F})) = \Phi_{r_0}(\pi_1(F_0))$
	
	Let $\iota_\ast \colon \pi_1(\overline{F}) \to \pi_1(F_0)$ be the homomorphism induced by the inclusion  $\iota \colon \overline{F} \to F_0$.
	Since $\overline{F}$ can be obtained from $F_0$ by deleting the interiors of some disks, every loop on $F_0$ can be homotoped to be a loop contained in $\overline{F}$. Thus $\iota_\ast$ is surjective. 
	
	We also claim that $\Phi_{r_0} = \Phi_{\overline{r}} \circ \iota_\ast$. To see this,  recall that the permutation representations are found by lifting loops to paths in the covers. Since $\overline{F} \subset F_0$, $\overline{S}\subseteq S_0$, and $\overline{r}$ is the restriction of $r_0$ to $\overline{S}$, every loop in $\overline{F}$ will lift to the same paths in $S_0$  whether or not you first include $\overline{F}$ into $F_0$. Hence, $\Phi_{r_0} (\gamma) = \Phi_{\overline{r}} \circ \iota_\ast(\gamma)$ for each $\gamma \in \pi_1(\overline{F})$.
	
	Because $\iota_\ast$ is surjective and $\Phi_{r_0} = \Phi_{\overline{r}} \circ \iota_\ast$, we must have that $\Phi_{\overline{r}}(\pi_1(\overline{F})) = \Phi_{r_0}(\pi_1(F_0))$ as desired.
\end{proof}

\subsection{Orbifold homeomorphisms acting on covers}\label{sec:orbifold_action}

Let $p \colon S \to \orb$ be a finite degree $D$ branched cover of an orbifold $\orb$ by a surface $S$ where the branch points are exactly the orbifold points of $\orb$. Let $\overline{O}$ be the surface obtained from $O$ by deleting a small neighborhood of each orbifold point as in Section \ref{subsec:permutations}.

Let $\orb$ and $\orbprm$ be two orbifolds with the same signature. We say two degree $D$ branched covers $p\colon S \to \orb$ and $p' \colon S\to \orbprm$ are \emph{signature equivalent} if there exists homeomorphisms $h \colon S \to S$ and $f \colon \orb \to \orbprm$ so that $f \circ p = p' \circ h$. Note, this is equivalent to saying that  $f \circ p$ is covering space equivalent to $p'$.

Now consider covers of orbifolds with signature $(0;r_1,r_2,r_3,r_4)$ where the orbifold point $x_i$ has order $r_i$. Fix generators $\gamma_1,\gamma_2,\gamma_3,\gamma_4$ for $\pi_1(\overline{O})$ as shown in in Figure \ref{fig:4-holed_sphere}. We will adopt the convention that $\delta_i$ is the boundary of a small neighborhood of the order $r_i$ orbifold point $x_i$. The permutation representation of a degree $D$ branched cover $p$  are determined by tuples $[s_1,s_2,s_3,s_4]$ of element of $\sym(D)$ where $s_i = \Phi_p(\gamma_i)$; see Section \ref{subsec:tuples}. If $f \colon \orb \to \orbprm$ is an orbifold homeomorphism that preserves the basepoint of $\pi_1(\overline{O})$, then $\Phi_{f \circ p}$ and  $\Phi_{p}$ are related by $\Phi_{f\circ p} =  \Phi_{p} \circ f^{-1}_\ast$ where $f_\ast$ is the automorphisms of $\pi_1(\overline{O})$ induced by $f$. This induces an action of  $f$ on the $\sym(D)$-tuples that is given by $$f \cdot [s_1,s_2,s_3,s_4] =  [ \Phi_p(f_\ast^{-1}(\gamma_1)), \Phi_p(f_\ast^{-1}(\gamma_2)), \Phi_p(f_\ast^{-1}(\gamma_3)), \Phi_p(f_\ast^{-1}(\gamma_4)) ].$$

If $f$ does not preserve the basepoint of $\pi_1(\overline{O})$, there is an orbifold homeomorphism $g$ of $\orbprm$ so that $g \circ f$ fixes the base point of $\pi_1(\overline{O})$ and $g$ is isotopic to the identity. If $g'$ is a different choice of homeomorphism of $\orbprm$ so that $g'$ is isotopic to the identity and  $g'\circ f$ fixes the basepoint of $\pi_1(\overline{O})$, then the action of  $g\circ f$ and $g'\circ f$ on $\pi_1(\overline{O})$ will differ by conjugation in $\pi_1(\overline{O})$. Hence there is an action of $\MCG\orb$ on the set of conjugacy classes of the $\sym(D)$-tuples by having each mapping class act by any choice of homeomorphism that fixes the basepoint.  

The above discussion applies equally well for orbifold with any signature. We will also apply it in the case of orbifolds with signature $(1;r)$ with the tuples discussed in Section \ref{subsec:tuples}.

The next lemma describes how we can classify when coverings are signature equivalent by using the $\MCG\orb$-orbits of their associated tuples.

\begin{lemma}\label{lemma:MCG(0)-orbits}
	Let $\orb$ and $\orbprm$ be orbifolds whose  signatures are either both $(0;r_1,r_2,r_3,r_4)$ or both $(1;r)$. Suppose $p\colon S \to \orb$ and $p'\colon S \to \orbprm$ are two degree $D$ branched covering maps. If $\bt$ and $\bt'$ are the $\sym(D)$-tuples for the permutation representations of $p$ and $p'$ respectively,  then $p$ and $p'$ are signature equivalent if and only if the $\MCG\orb$-orbit of the conjugacy class of $\bt$ contains the conjugacy class of $\bt'$.
\end{lemma}

\begin{proof}
	Start by assuming that $p$ and $p'$ are signature equivalent. There then exists homeomorphisms $h \colon S \to S$ and $f \colon \orb \to \orbprm$ so that $p' \circ h = f \circ p$. Since the adjustment to make $f$ fix the basepoint of $\pi_1(\overline{O})$ is isotopic to the identity, it lifts to $S$ and we can assume that $f$ fixes the basepoint of  $\pi_1(\overline{O})$ without changing the isotopy class of either $h$ or $f$.  With this adjustment, $p$ and $p'$ being signature equivalent is equivalent to $\Phi_{f\circ p}$ and $\Phi_p$ being conjugate in $\sym(D)$. Since different choices for the adjustment only change  $\Phi_{f\circ p}$ be conjugation, we have  $f \cdot \bt$ is conjugate to $\bt'$ for any choice of adjustment.
	
	Now assume there is an orbifold homeomorphism $f$ that fixes the basepoint of $\pi_1(\overline{O})$  and acts on $\bt$ so that $f \cdot \bt$ is conjugate to $\bt'$. Thus $f\circ p$ is covering space equivalent to $p'$, which implies $p$ and $p'$ are signature equivalent.
\end{proof}

Because of Lemma \ref{lemma:MCG(0)-orbits}, if the conjugacy classes of the $\sym(D)$-tuples  $\bt$ and $\bt'$ are  in the same $\MCG\orb$-orbit, we say $\bt$ and $\bt'$ are 
\emph{signature equivalent}. An important invariant of signature equivalent tuples is the isomorphism class of the subgroup of $\sym(D)$ that they generate.

\begin{corollary}
	Suppose $\bt$ and $\bt'$ are signature equivalent $\sym(D)$-tuples for signatures $(0;r_1,r_2,r_3,r_4)$ or $(1;r)$. The subgroup of $\sym(D)$ that is generated by the element of $\bt$ is isomorphic to the subgroup generated by $\bt'$.
\end{corollary}

\begin{figure}[ht]
	\centering
    \def\svgwidth{3in}
\begingroup%
  \makeatletter%
  \providecommand\color[2][]{%
    \errmessage{(Inkscape) Color is used for the text in Inkscape, but the package 'color.sty' is not loaded}%
    \renewcommand\color[2][]{}%
  }%
  \providecommand\transparent[1]{%
    \errmessage{(Inkscape) Transparency is used (non-zero) for the text in Inkscape, but the package 'transparent.sty' is not loaded}%
    \renewcommand\transparent[1]{}%
  }%
  \providecommand\rotatebox[2]{#2}%
  \newcommand*\fsize{\dimexpr\f@size pt\relax}%
  \newcommand*\lineheight[1]{\fontsize{\fsize}{#1\fsize}\selectfont}%
  \ifx\svgwidth\undefined%
    \setlength{\unitlength}{148.74064514bp}%
    \ifx\svgscale\undefined%
      \relax%
    \else%
      \setlength{\unitlength}{\unitlength * \real{\svgscale}}%
    \fi%
  \else%
    \setlength{\unitlength}{\svgwidth}%
  \fi%
  \global\let\svgwidth\undefined%
  \global\let\svgscale\undefined%
  \makeatother%
  \begin{picture}(1,0.30999378)%
    \lineheight{1}%
    \setlength\tabcolsep{0pt}%
    \put(0.46412917,0.16140757){\makebox(0,0)[lt]{\lineheight{1.25}\smash{\begin{tabular}[t]{l}$H_{1,2}$\end{tabular}}}}%
    \put(0,0){\includegraphics[width=\unitlength,page=1]{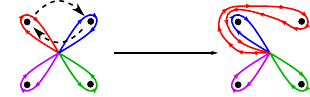}}%
    \put(-0.00549716,0.23148226){\color[rgb]{1,0,0}\makebox(0,0)[lt]{\lineheight{1.25}\smash{\begin{tabular}[t]{l}$\gamma_1$\end{tabular}}}}%
    \put(0.32441206,0.22524163){\color[rgb]{0,0,1}\makebox(0,0)[lt]{\lineheight{1.25}\smash{\begin{tabular}[t]{l}$\gamma_2$\end{tabular}}}}%
    \put(-0.00669684,0.02120091){\color[rgb]{0.78431373,0,1}\makebox(0,0)[lt]{\lineheight{1.25}\smash{\begin{tabular}[t]{l}$\gamma_4$\end{tabular}}}}%
    \put(0.332926,0.01902692){\color[rgb]{0,0.70588235,0}\makebox(0,0)[lt]{\lineheight{1.25}\smash{\begin{tabular}[t]{l}$\gamma_3$\end{tabular}}}}%
  \end{picture}%
\endgroup%

    \def\svgwidth{3in}
\begingroup%
  \makeatletter%
  \providecommand\color[2][]{%
    \errmessage{(Inkscape) Color is used for the text in Inkscape, but the package 'color.sty' is not loaded}%
    \renewcommand\color[2][]{}%
  }%
  \providecommand\transparent[1]{%
    \errmessage{(Inkscape) Transparency is used (non-zero) for the text in Inkscape, but the package 'transparent.sty' is not loaded}%
    \renewcommand\transparent[1]{}%
  }%
  \providecommand\rotatebox[2]{#2}%
  \newcommand*\fsize{\dimexpr\f@size pt\relax}%
  \newcommand*\lineheight[1]{\fontsize{\fsize}{#1\fsize}\selectfont}%
  \ifx\svgwidth\undefined%
    \setlength{\unitlength}{156.83797602bp}%
    \ifx\svgscale\undefined%
      \relax%
    \else%
      \setlength{\unitlength}{\unitlength * \real{\svgscale}}%
    \fi%
  \else%
    \setlength{\unitlength}{\svgwidth}%
  \fi%
  \global\let\svgwidth\undefined%
  \global\let\svgscale\undefined%
  \makeatother%
  \begin{picture}(1,0.33128971)%
    \lineheight{1}%
    \setlength\tabcolsep{0pt}%
    \put(0.44605243,0.19059556){\makebox(0,0)[lt]{\lineheight{1.25}\smash{\begin{tabular}[t]{l}$H_{1,3}$\end{tabular}}}}%
    \put(0,0){\includegraphics[width=\unitlength,page=1]{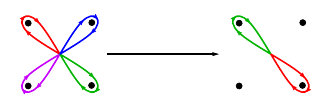}}%
    \put(-0.0063511,0.23626978){\color[rgb]{1,0,0}\makebox(0,0)[lt]{\lineheight{1.25}\smash{\begin{tabular}[t]{l}$\gamma_1$\end{tabular}}}}%
    \put(0.29696134,0.23035134){\color[rgb]{0,0,1}\makebox(0,0)[lt]{\lineheight{1.25}\smash{\begin{tabular}[t]{l}$\gamma_2$\end{tabular}}}}%
    \put(0.00207518,0.03684497){\color[rgb]{0.78431373,0,1}\makebox(0,0)[lt]{\lineheight{1.25}\smash{\begin{tabular}[t]{l}$\gamma_4$\end{tabular}}}}%
    \put(0.32354021,0.05166886){\color[rgb]{0,0.70588235,0}\makebox(0,0)[lt]{\lineheight{1.25}\smash{\begin{tabular}[t]{l}$\gamma_3$\end{tabular}}}}%
    \put(0,0){\includegraphics[width=\unitlength,page=2]{4_loops_across_twist_svg-tex.pdf}}%
  \end{picture}%
\endgroup%

	
	\caption{Examples of the homeomorphisms $H_{i,j}$. All other $H_{i,j}$ can be obtained from the shown examples by rotating or using the fact that $H_{j,i} = H^{-1}_{i,j}$.} 
	\label{fig:half_twists}
\end{figure}

To understand the $\MCG\orb$-orbits of the $\sym(D)$-tuples, we focus on some simple  ``braiding'' homeomorphisms of  $(0;r_1,r_2,r_3,r_4)$ orbifolds that will generate $\MCG\orb$. For distinct $i,j \in \{1,2,3,4\}$, let $H_{i,j}$ denote the the surface homeomorphism of $O$ shown in Figure \ref{fig:half_twists} that exchanges the orbifold points $x_i$ and $x_j$. Let $F_{i,j} = H^2_{i,j}$. Note $F_{i,j}$ is always an orbifold homeomorphism of $\orb$, while $H_{i,j}$ will be an orbifold homeomorphism if and only if $\order(x_i) = \order(x_j)$. The mapping class group $\MCG\orb$ is well known to be generated by the set of $H_{i,j}$ and $F_{i,j}$ that it contains; see for example \cite{FarbMarg}.

Figure \ref{fig:half_twists} shows the action of $H_{i,j}$ on the loops $\gamma_1,\gamma_2,\gamma_3,\gamma_4$ that generate $\pi_1(\overline{O})$. The action of $H_{i,j}$ and $F_{i,j}$ on the $\gamma_i$ are summarized in the following tables where all indicies are taken modulo 4. For the remainder of the paper, we will use an overline on a group element to indicate the inverse of that group element. That is, $\overline{\gamma} = \gamma^{-1}$.

\begin{multicols}{2}
	
	\begin{tabular}{c |c |c}
		& $H_{i,i+1}$& $F_{i,i+1} = H^2_{i,i+1}$  \\ \hline
		$\gamma_i$ & ${\gamma}_i {\gamma}_{i+1} \overline{\gamma}_i$ & $ ({\gamma}_i {\gamma}_{i+1}) {\gamma}_i (\overline{\gamma}_{i+1} \overline{\gamma}_i) $ \\
		
		$\gamma_{i+1}$ &  $\gamma_i $ &  ${\gamma}_i {\gamma}_{i+1} \overline{\gamma}_i$\\
		
	\end{tabular}
	
	\vspace{.5cm}

	\begin{tabular}{c |c |c}
		& $H_{i,i+2}$& $F_{i,i+2} = H^2_{i,i+2}$  \\ \hline
		$\gamma_i$  & $\gamma_{i+2}$ &  $\gamma_i$\\
		
		$\gamma_{i+1}$ & $\overline{\gamma}_{i+2} \gamma_{i+1} {\gamma}_{i+2} $ & $ (\overline{\gamma}_i \overline{\gamma}_{i+2}) \gamma_{i+1} (\gamma_{i+2} \gamma_{i})$
		\\
		
		$\gamma_{i+2}$ & $\gamma_i$ & $\gamma_{i+2}$\\
		
		$\gamma_{i+3}$ & $\overline{\gamma}_{i} \gamma_{i+3} \gamma_i$ & $(\overline{\gamma}_{i+2} \overline{\gamma}_i) \gamma_{i+3} (\gamma_{i} \gamma_{i+2}) $
		\\
		
	\end{tabular}
	
	\begin{tabular}{c |c | c}
		& $H_{i+1,i}$& $F_{i+1,i} = H^2_{i+1,i}$  \\\hline
		
		$\gamma_i$ & $\gamma_{i+1}$ & $\overline{\gamma}_{i+1} \gamma_i  {\gamma}_{i+1}$ \\
		
		$\gamma_{i+1}$ & $\overline{\gamma}_{i+1} \gamma_i {\gamma}_{i+1}$ &  $ (\overline{\gamma}_{i+1}\overline{\gamma}_i) \gamma_{i+1} ({\gamma}_i {\gamma}_{i+1})$\\
	\end{tabular}
	\vspace{.5cm}
	
	\begin{tabular}{c |c |c}
		& $H_{i+2,i}$& $F_{i+2,i} = H^2_{i+2,i}$  \\ \hline

		$\gamma_i$  & $\gamma_{i+2}$ &  $\gamma_i$\\
		
		$\gamma_{i+1}$ & $\gamma_{i} {\gamma}_{i+1}  \overline{\gamma}_{i}$ & $(\gamma_{i+2} \gamma_i) \gamma_{i+1} (\overline{\gamma}_i \overline{\gamma}_{i+2})$
		\\
		
		$\gamma_{i+2}$ & $\gamma_i$ & $\gamma_{i+2}$\\
		
		$\gamma_{i+3}$ & $ \gamma_{i+2} {\gamma}_{i+3} \overline{\gamma}_{i+2}$ & $(\gamma_i \gamma_{i+2}) \gamma_{i+3} (\overline{\gamma}_{i+2} \overline{\gamma}_i) $
		\\
	\end{tabular}
\end{multicols}
Since $\Phi_{H_{i,j} \circ p} = \Phi_p \circ H_{i,j}^{-1}$ and $H_{i,j}^{-1} = H_{j,i}$, the following tables describe the action of $H_{i,j}$ and $F_{i,j}$ on a tuple $[s_1,s_2,s_3,s_4]$ for a branched covering $p$.

\begin{multicols}{2}

	\begin{tabular}{c |c |c}
		& $H_{i,i+1}$& $F_{i,i+1} = H^2_{i,i+1}$  \\ \hline
		$s_i$ & $s_{i+1}$ & $ \overline{s}_{i+1} s_i s_{i+1}$ \\
		
		$s_{i+1}$ & $\overline{s}_{i+1} s_i s_{i+1}$ &  $ (\overline{s}_{i+1}\overline{s}_i) s_{i+1} ({s}_i s_{i+1})$\\
		
	\end{tabular}
	
	\vspace{.5cm}

	\begin{tabular}{c |c |c}
		& $H_{i,i+2}$& $F_{i,i+2} = H^2_{i,i+2}$  \\ \hline
		$s_i$  & $s_{i+2}$ &  $s_i$\\
		
		$s_{i+1}$ & $s_{i} {s}_{i+1}  \overline{s}_{i}$ & $(s_{i+2} s_i) s_{i+1} (\overline{s}_i \overline{s}_{i+2})$\\
		
		$s_{i+2}$ & $s_i$ & $s_{i+2}$\\
		
		$s_{i+3}$ & $ s_{i+2} {s}_{i+3} \overline{s}_{i+2}$ & $(s_i s_{i+2}) s_{i+3} (\overline{s}_{i+2} \overline{s}_i) $\\
	\end{tabular}
	
	\columnbreak
	
	\begin{tabular}{c |c | c}
		& $H_{i+1,i}$& $F_{i+1,i} = H^2_{i+1,i}$  \\\hline
		$s_i$ & $s_i {s}_{i+1} \overline{s}_i$ & $ (s_i {s}_{i+1}) {s}_i (\overline{s}_{i+1} \overline{s}_i) $ \\
		
		$s_{i+1}$ &  $s_i $ &  $s_i {s}_{i+1} \overline{s}_i$\\
	\end{tabular}
	
	\vspace{.5cm}

	\begin{tabular}{c |c |c}
		& $H_{i+2,i}$& $F_{i+2,i} = H^2_{i+2,i}$  \\ \hline
		$s_i$  & $s_{i+2}$ &  $s_i$\\
		
		$s_{i+1}$ & $\overline{s}_{i+2} s_{i+1} {s}_{i+2} $ & $ \overline{s}_i \overline{s}_{i+2} s_{i+1} s_{i+2} s_{i+1}$\\
		
		$s_{i+2}$ & $s_i$ & $s_{i+2}$\\
		
		$s_{i+3}$ & $\overline{s}_{i} s_{i+3} s_i$ & $\overline{s}_{i+2} \overline{s}_i s_{i+3} s_{i} s_{i+2} $\\
	\end{tabular}
\end{multicols}

\section{Flexible cone metrics and covers}\label{sec:flexible_cone_metrics}
Hersonsky and Paulin showed that a metric $\varphi \in \hyp(S)$ is determined by its associated \emph{Louisville geodesic current} \cite{HP_marked_length}. Erlandsson, Leininger, and Sadanand showed that in all but very specific cases only the \emph{support}, $\G_\varphi$, of the Louisville current is need to determine the metric \cite{ELS_hyp_cone_metrics}. We say a metric $\varphi\in\hyp(S)$ is \emph{flexible} if there is a metric $\varphi' \in \hyp(S) - \{\varphi\}$ so that $\G_\varphi = \G_{\varphi'}$. If $\varphi$ is not flexible, then we say it is \emph{rigid}. When $\varphi$ is flexible, the \emph{flexibility class} of $\varphi$ is the set of metrics $\varphi' \in \hyp(S)$ with $\G_{\varphi'}=\G_\varphi$. We use $\mathcal{F}_\varphi$ to denote the flexibility class of $\varphi$.

Erlandsson, Leininger, and Sadanand understand flexible metrics using specific branched covers of orbifolds. We say a map $p \colon (S,\varphi) \to \orb$ is a \emph{flexible cover} if
\begin{enumerate}
	\item $\orb$ is a non-triangular\footnote{A triangular orbifold is one with signature $(0;r_1,r_2,r_3)$.} hyperbolic orbifold,
	\item  $p$ is a finite degree, locally isometric branched cover, 
	\item each cone point in $(S,\varphi)$ is sent to an even order orbifold point of $\orb$.
\end{enumerate}

When a metric admits a flexible cover, varying the metric on the covered orbifold produces distinct metrics with the same Louisville current. 

\begin{theorem}\label{thm:flexible_covers_give_flexible_metrics}
	Let $S$ be a closed, orientable surface with genus at least 2 and let $\varphi \in \hyp(S)$. If $p\colon (S,\varphi) \to \orb$ is a flexible covering map, then for every orbifold homeomorphism $f \colon \orb \to \orbprm$, the metric $\varphi'$ obtained by pulling back the metric on $\orbprm$ under $f \circ p$ has $\mathcal{G}_\varphi = \mathcal{G}_{\varphi'}$. Moreover, $\varphi \neq \varphi'$ if and only if $f$ is not orbifold isotopic to an isometry.
\end{theorem}

Erlandsson, Leininger, and Sadanand show that the entire flexibility class of a flexible metric can be described by constructing a canonical flexible covering and then deforming the metrics on the this orbifold. We recall this construction, and demonstrate some of its consequences.

\begin{definition}[Holonomy cover]
	Given a flexible metric $\varphi \in \hyp(S)$, let $(\hat{S},\hat{\varphi})$ be the metric completion of the universal cover of the surface obtained by puncturing $(S,\varphi)$ at $\cone(\varphi)$. There is a $\pi_1(\dot{S})$-equivariant covering map $\hat{j} \colon (\hat{S},\hat{\varphi}) \to (S,\varphi)$, a developing map $D \colon (\hat{S},\hat{\varphi}) \to \HP$, and holonomy representation $\rho\colon \pi_1(\dot{S}) \to \PSL(2,\R)$ as described in Subsection \ref{subsec:punctures}.  For each $\hat{x} \in \hat{j}^{-1}(\cone(\varphi))$, let $x = D(\hat{x})$ and $\tau_x \in \PSL(2,\R)$ be the order two rotation around $x$.  Define $\Lambda$ to be the subgroup of $\PSL(2,\R)$ generated by $\rho(\pi_1(\dot{S}))$ and the set $\{\tau_{x} : x \in \hat{j}^{-1}(\cone(\varphi))\}$.  Erlandsson, Leininger, and Sadanand prove that when $\varphi$ is flexible, $\Lambda$ is a discrete subgroup. Thus $\orb = \HP / \Lambda$ is a hyperbolic orbifold and there is a commuting diagram
	
	\begin{center}
		\begin{tikzcd}
			(\hat{S},\hat{\varphi}) \arrow[r,"\hat{j}"] \arrow[d, "D"] & (S,\varphi) \arrow[d, "p"] \\
			\HP \arrow[r] & \orb
		\end{tikzcd}
	\end{center}
	
	\noindent where $p$ is a flexible covering.
	
	Since this construction depends only on the choice of developing map $D$, we call $\Lambda$ the \emph{holonomy lattice for $D$},  $\orb = \HP /\Lambda$ the \emph{holonomy orbifolds for $D$}, and $p \colon (S,\varphi) \to \orb$ the \emph{holonomy cover for $D$}.
	Recall, if $D'$ is another developing map for $\varphi$, then there is $f \in \PSL(2,\R)$ so that $D' = f\circ D$. Hence $f^{-1} \Lambda f$ will be the holonomy lattice for $D'$ and $f$ will descend to a isometry between the holonomy orbifolds for $D$ and $D'$. This isometry commutes with the two holonomy covers. Hence, when the specific developing map is not relevant, we will refer simply to a holonomy lattice/orbifold/cover.
\end{definition}

\begin{theorem}[\cite{ELS_hyp_cone_metrics}]\label{thm:flexible_metrics_give_flex_covers}
	Let $S$ be a closed, orientable surface with genus at least 2. Suppose $\varphi \in \hyp(S)$ is a flexible metric and $\varphi' \in \mathcal{F}_\varphi - \{\varphi\}$. 
	\begin{enumerate}
		\item The holonomy covers $p\colon (S,\varphi) \to \orb$ and $p'\colon (S,\varphi') \to \orbprm$  for $\varphi$ and $\varphi'$ are flexible covers.
		\item There exists  homeomorphisms $h \colon S \to S$ and  $f \colon \orb \to \orbprm$ so that $p' \circ h = f \circ p$, $h$  is isotopic to the identity relative $\cone(\varphi)$, and $\varphi'$ is the pull-back of the metric on $\orbprm$ under $f \circ p$.
	\end{enumerate}

\end{theorem}

Theorem \ref{thm:flexible_metrics_give_flex_covers} implies that the holonomy orbifold is the maximal flexible covering possible for a given metric.

\begin{corollary}\label{cor:maximal_flexible_cover}
	Suppose $\varphi \in \hyp(S)$ is flexible. Let  $ p_0 \colon (S,\varphi) \to \orbnot$ be a holonomy cover for $\varphi$.  If $p \colon (S,\varphi) \to \orb$ is a flexible covering, then there exists a locally isometric branch cover $q \colon \orbnot \to \orb$ so that $p = q \circ p_0$.
\end{corollary}

\begin{proof}
	Let $\Lambda_0$ be the holonomy lattice for a fix choice of developing map $D$ for $(\hat{S},\hat{\varphi})$. Let $\orbnot = \HP / \Lambda_0$ be the holonomy orbifold and $p_0 \colon (S,\varphi) \to \orbnot$ be the holonomy covering. Let $\dot{S}$ be the topological surface obtained by puncturing $S$ at $\cone(\varphi)$.

	Let $\Lambda$ be the discrete subgroup of $\PSL(2,\R)$ so that $\orb = \HP / \Lambda$. Since there is a locally isometric branched cover $p \colon (S,\varphi) \to \orb$, there is a local isometry from $(\hat{S},\hat{\varphi}) \to \HP$ that makes the following diagram commutes:
	 
	\begin{center}
		\begin{tikzcd}
			(\hat{S},\hat{\varphi}) \arrow[r],\arrow[d] & (S,\varphi)\arrow[d,"p"]\\ 
			\HP \arrow[r] & \orb = \HP / \Lambda
		\end{tikzcd}
	\end{center}
	\noindent By conjugating $\Lambda$ in $\PSL(2,\R)$, we can take the map $(\hat{S},\hat{\varphi}) \to \HP$ to be our fixed developing map $D$. 
	
	Now, $p \colon (S,\varphi) \to \orb$ gives a representation $\rho' \colon \pi_1(\dot{S}) \to \PSL(2,\R)$ that is equivariant with the action of $\pi_1(\dot{S})$ on $\hat{S}$. Hence for all $\gamma \in \pi_1(\dot{S})$, $$\rho'(\gamma)  \circ D = D \circ \gamma  = \rho(\gamma) \circ D.$$ Thus $\rho' = \rho$. This implies $\rho(\pi_1(\dot{S}))$ is a subgroup of $\Lambda$. 
	
	Let $ \hat{x} \in \hat{j}^{-1}(\cone(\varphi))$ and $x = D(\hat{x})$. Since $p$ is flexible, $p(\hat{j}(\hat{x}))$ must be an even order orbifold point of $\orb$. Thus, $\tau_x$---the order two rotation around $x$ in $\HP$---must be contained in $\Lambda$. Specially $\tau_x \in \Lambda$ will be a power of a loop that goes around the puncture at $p(\hat{j}(\hat{x}))$. Thus, $\Lambda_0 < \Lambda$ as $\Lambda$ contains a generating set for $\Lambda_0$.     
	By covering space theory, $\Lambda_0 < \Lambda$ implies there is a locally isometric branched covering map $q \colon \orbnot \to \orb$ so that $p = q \circ p_0$. 
\end{proof}

When two flexible metrics are in the same $\MCG(S)$-orbit, their holonomy covers differ by any isometry of the covered orbifold.

\begin{lemma}\label{lem:orbits_of_flex_covers}
	Suppose $\varphi_1,\varphi_2 \in \hyp(S)$ are flexible metrics and  there is a homeomorphism $h \colon S \to S$ so that $h(\varphi_1) = \varphi_2$. Let $p_i \colon (S,\varphi_i) \to \orbi$ be a holonomy cover for $\varphi_i$. There exists an isometry $F \colon \orbone \to \orbtwo$ so that $p_2 \circ h = F \circ p_1$.

\end{lemma}

\begin{proof}
	Since $h(\varphi_1) = \varphi_2$, $\varphi_1$ and $\varphi_2$ have the same number of cone points. Thus, puncturing $S$ at $\cone(\varphi_1)$ or $\cone(\varphi_2)$ will produce the same topological surface $\dot{S}$.  We can also adjust $h$ by a homeomorphism that is isotopic to the identity so that $h$ will fix the base point of $\pi_1(S)$. With this modifications, $h$ induces an automorphism $h_\ast \colon \pi_1(\dot{S}) \to \pi_1(\dot{S})$ that is well defined up to conjugation.
	
	The homeomorphism $h$ lifts to a $h_\ast$--equivariant homeomorphism $\hat{h} \colon (\hat {S}, \hat{\varphi_1}) \to (\hat {S}, \hat{\varphi_2})$.  Since $h(\varphi_1) = \varphi_2$, there exists a pair of developing maps $D_i \colon (\hat {S}, \hat{\varphi_i}) \to \HP$ so that $D_1 = D_2 \circ \widehat{h}$. For each of these developing maps, there is a representation $\rho_i \colon \pi_1(\dot{S}) \to \PSL(2,\R)$ with the property $$D_i \circ \gamma  = \rho_i(\gamma) \circ D_i$$ for each $ \gamma \in \pi_1(\dot{S})$.
	
	The following calculation verifies that $\rho_1 = \rho_2 \circ h_\ast$
	
	\[ \rho_1(\gamma) \circ D_1 = D_1 \circ \gamma = D_2 \circ  \widehat{h} \circ \gamma  = D_2 \circ h_\ast(\gamma) \circ  \widehat{h} = \rho_2( h_\ast(\gamma)) \circ  (D_2 \circ \widehat{h}) = \rho_2(h_\ast(\gamma)) \circ D_1. \]
	
	Since $h_\ast$ is an automorphism of $\pi_1(\dot{S})$, the images of $\rho_1$ and $\rho_2$ in $\PSL(2,\R)$ are equal. Hence, the holonomy lattices for $D_1$ and $D_2$ are equal. Thus, there is a single holonomy orbifold $\orbone$ for $D_1$ and $D_2$. In particular, we have the following commuting diagram:
	
	\begin{center}
		\begin{tikzcd}[sep=small]
			(\hat{S},\hat{\varphi}_1)  \arrow[dr, "D_1" below] \arrow[rr,"\hat{h}"] & &(\hat{S},\hat{\varphi}_2) \arrow[ld,"D_2" ] \\
			& \HP \arrow[d] &\\
			& \orbone & 
		\end{tikzcd}
	\end{center}
	
	Now the holonomy orbifold $\orbtwo$ and the flexible covering $p_2 \colon (S,\varphi_2) \to \orbtwo$ are determined by some choice of developing map $D'_2$ for $\varphi_2$. Thus there is a isometry $\widetilde{F} \colon \HP \to \HP$ so that $D_2' = \widetilde{F} \circ D_2$.  This isometry descend to an isometry $F \colon \orbone \to \orbtwo$ so that $p_2=F \circ p$. Hence $$p_2 \circ h = p_1 \circ F.$$
\end{proof}
Combining these results, we see that the $\MCG(S)$-orbits of flexibility classes in $\hyp(S)$ correspond to the holonomy covering maps modulo homeomorphisms of both the surface and the covered orbifold.

\begin{corollary}\label{cor:flex_class_iff_MCG-equiv}
	Let $\varphi_1,\varphi_2 \in \hyp(S)$ be two flexible metrics. Let $p_i\colon (S,\varphi_i) \to \orbi$ be a holonomy cover for $\varphi_i$. 
	The flexibility class of $\varphi_1$ is in the same $\MCG(S)$-orbit as the flexibility class of $\varphi_2$ if  and only if $p_1$ and $p_2$ are signature equivalent.
\end{corollary}

\begin{proof}
	Let $\mathcal{F}_i$ be the flexibility class  $\mathcal{F}_{\varphi_i}$. 
	
	First assume there exists a homeomorphism $h \colon S \to S$ so that $h(\mathcal{F}_1) = \mathcal{F}_2$. Hence $\varphi'_2 = h(\varphi_1)$ is in $\mathcal{F}_2$. By Theorem \ref{thm:flexible_covers_give_flexible_metrics}, there exists an orbifold $(O_2,\psi_2')$ and a homeomorphism $f_2 \colon \orbtwo \to (O_2,\psi_2')$ so that $\varphi'_2$ is the pull-back of $\psi_2$ under $f_2\circ p_2$. Now Lemma \ref{lem:orbits_of_flex_covers} say there exists an isometry $F \colon \orbone \to (O_2,\psi_2')$ so that $F \circ p_1 = (f_2\circ p_2) \circ h$. The homeomorphism $$f = f_2^{-1} \circ F \colon \orbone \to \orbtwo,$$ therefore has $f \circ p_1 = p_2 \circ h$.
	
	Now assume $p_1$ and $p_2$ are signature equivalent. That is, there exists homeomorphisms $h \colon S \to S$ and $f \colon \orbone \to \orbtwo$  so that $p_2 \circ h = f \circ p_1$.  Let $\varphi'_1$ be the pull back of the metric on $\orbtwo$ under $f \circ p_1$.  Theorem \ref{thm:flexible_covers_give_flexible_metrics} says that $\varphi'_1 \in \mathcal{F}_1$. Now the metric $h(\varphi'_1)$ is the pull-back of the metric on $\orbtwo$ under the map $f \circ p_1 \circ h^{-1}$. Since $p_2 = f\circ p_1 \circ h^{-1}$ by assumption, $\varphi_2 = h(\varphi'_1)$. Hence $\mathcal{F}_2 = h(\mathcal{F}_1)$.
\end{proof}

\subsection{Some lemmas on flexible covers}\label{section:basic_lemmas}
Fix a surface $S$ and a metric $\varphi \in \hyp(S)$. Let $p\colon (S,\varphi) \to \orb$ be  a flexible covering map with finite degree $D$.  We establish a few basic facts about $p$. In the next section, we will use these facts to limit the orbifolds that can have a flexible covering by a genus 2 surface.

\begin{lemma}\label{lem:order(x)<D}
	For every orbifold point $x\in \orb$, $\order(x) \leq D$.
\end{lemma}
\begin{proof}
	
	Let $x \in \orb$ be an orbifold point and $y \in p^{-1}(x)$. By Equation \ref{eq:cover_degree}, 
	\[ \Theta_y =\deg(y) \Theta_x  \implies  2\pi  \leq  \deg(y)\frac{2\pi}{\order(x)}   
	\implies   \order(x) \leq  \deg(y). 
	\]  Since $\deg(y) \leq D$, we have $\order(x) \leq D$ as desired.
\end{proof}

\begin{lemma}\label{lem:odd|D}
	If $x\in \orb$ is an orbifold point and $\order(x)$ is odd, then $\order(x)$ divides $D$.
\end{lemma}
\begin{proof}        
	Since $p$ is flexible, every preimage of an odd order orbifold point of $\orb$ is a regular point of $(S,\varphi)$. Thus, for all odd order orbifold points $x \in \orb$ and $ y\in p^{-1}(x)$, we have $\Theta_y = 2\pi$.
	By Equation \ref{eq:cover_degree}, 
	\[\Theta_x = \frac{\Theta_y}{\deg(y)} \implies
	\frac{2\pi}{\order(x)}  = \frac{2\pi}{\deg(y)} \implies
	\order(x) = \deg(y).       \]
	By Equation \ref{eq:sum_local_degree}, we get
	\[D=\sum_{y\in p^{-1}(x)}\deg(y) =\sum_{y\in p^{-1}(x)}\order(x) = |p^{-1}(x)| \cdot \order(x).\]
	Thus, $\order(x)$ divides $D$.                    
\end{proof}

\begin{lemma}\label{regular preimages}
	Suppose $x\in \orb$ is an even order orbifold point. If every $y \in p^{-1}(x)$ is a regular point of $(S,\varphi)$, then $D$ must be even and  $\deg(y) = \order(x)$ for all $y \in p^{-1}(x)$. 
\end{lemma}
\begin{proof}
	Let $x\in \orb$ be an orbifold point with $\order(x) = 2k$ for some integer $k$. Assume every point in $p^{-1}(x)$ is a regular point of $(S,\varphi)$. By Equation \ref{eq:cover_degree}, for all  $y \in p^{-1}(x)$, $$\deg(y) = \frac{\Theta_x}{\Theta_y} = \frac{2\pi/2k}{2\pi} = 2k = \order(x).$$
	By Equation \ref{eq:sum_local_degree}, we get 
	
	$$D= \sum_{y \in p^{-1}(x)}\deg(y)  \implies
	D= 2k |p^{-1}(x)|.$$
	Therefore, $D$ is even.
\end{proof}

\begin{lemma}\label{n_max}
	Suppose $S$ has genus $h$. Let $r$ be the order of the orbifold point of $\orb$ with the largest order among even order orbifold points of $\orb$. The number of cone points of $(S,\varphi)$ is bounded above by $n_{max} = -r\left( \frac{A(\psi)D - 4\pi h +4\pi}{2\pi}\right).$
\end{lemma}
\begin{remark}
	We will only apply Lemma \ref{n_max} when $h= 2$, so $n_{max} =  -r\left( \frac{A(\psi)D - 4\pi}{2\pi}\right)$.
\end{remark}

\begin{proof}
	Let $h$ be the genus of $S$ and $n$ be the number of cone points of $(S,\varphi)$. 
	To determine $n_{max}$, the upper bound for the number of cone points on $(S,\varphi)$, we  first determine the smallest possible angle of any cone point on $S$. 
	
	Let $y$ be a cone point on $(S,\varphi)$  and let $x \in \orb$ be $p(y) = x$. We  first determine a lower bound for $\Theta_y$.
	Since each cone point has angle strictly large than $2\pi$, we have $\deg(y)\Theta_{x} = \Theta_y > 2\pi.$  Thus, $\Theta_y$ is minimized when $\Theta_x$ is minimized. Since $x = p(y)$ must be an even order orbifold point of $\orb$, the smallest $\Theta_x$ can be is $2\pi/r$.  Hence, we have $\deg(y) \frac{2\pi}{r} > 2 \pi \implies \deg(y) \geq r+1,$ which implies $$\Theta_y \geq 2\pi \frac{r+1}{r}.$$
	
	Now let $\cone(\varphi)= (\Theta_1,\dots, \Theta_n)$. Plugging $\Theta_i \geq 2\pi(r+1)/r$ into the area of $\varphi$ in Equation \ref{Area Covering Theorem} produces
	\begin{align*}
		D &= \frac{-2\pi\left(2-2h-\sum_{i=1}^n\left(1-\Theta_i/2\pi\right)\right)}{A(\psi)}  \\
		D &\leq \frac{-2\pi\left(2-2h-\sum_{i=1}^{n}\left(1-\frac{r+1}{r}\right)\right)}{A(\psi)} \\
		D &\leq \frac{-2\pi\left(2-2h-n\left(-\frac{1}{r}\right)\right)}{A(\psi)} \\
		D &\leq \frac{4\pi h - 4\pi-2\pi n\frac{1}{r}}{A(\psi)} \\
		n & \leq -r\left( \frac{A(\psi)D - 4\pi h +4\pi}{2\pi}\right).\\
	\end{align*} 
	Hence,  $(S,\varphi)$ has at most $n_{max}=-r\left( \frac{A(\psi)D - 4\pi h +4\pi}{2\pi}\right)$ cone points.    
\end{proof}

\section{List of Possible Orbifolds} \label{Section: list of Orbifolds}
From now on, we will assume that our fixed surface $S$ has genus exactly 2. We continue to fix a metric $\varphi \in \hyp(S)$, and let $p\colon (S,\varphi) \to \orb$ be a flexible covering map with finite degree $D >1$. Let $(g;r_1,\dots,r_m)$ be the signature of $\orb$. By the end of this section, we will produce a list of all the possibilities for this signature.

First we determine the possibilities for the genus $g$ and number of orbifold points $m$ of $\orb$. 

\begin{lemma}\label{bounds on orbifold points}
	The orbifold $\orb$ has either $g=0$ and $m=4$, $g=0$ and $m=5$, or $g=1$ and $m=1$.
\end{lemma}
\begin{proof}     
	By Equation \ref{Area Covering Theorem}, we know that $D = A(\varphi)/A(\psi)$. A hyperbolic metric on a genus 2 surface with no cone points has area $4\pi$. As adding cone points with cone angle greater than $2\pi$ will always decrease the area, we have that
	\[ D \leq \frac{4\pi}{-2\pi(2-2g-\sum_{i=1}^m(1-\frac{1}{r_i}))} = \frac{2}{2g-2+ \sum_{i=1}^m(1-\frac{1}{r_i})}.    \]
	
	Since $D$ is a positive integer larger than 1, the denominator of the above expression must be between 0 and 1. That is
	
	$$0 <2g-2+\sum_{i=1}^m\left(1-\frac{1}{r_i}\right) < 1.$$
	This implies
	$$ 2-2g < \sum_{i=1}^m\left(1-\frac{1}{r_i}\right) < 3-2g.$$

	Now each $r_i \geq 2$ as each orbifold point has order at least 2. Thus we get, $$\sum_{i=1}^m(1-\frac{1}{r_i}) \geq \sum_{i=1}^m\left(1-\frac{1}{2}\right) = \frac{m}{2}.$$ This gives us the  bound $m<6-4g$.
	
	We can get the lower bound  $2-2g<m $ using the fact that  $\sum_{i=1}^m(1-\frac{1}{r_i})< \sum_{i=1}^m(1) = m$. Hence $$3-2g \leq m \leq 5-4g$$  as $m$ must be an integer.

	We now show $g = 0$ or $1$.  Since $A(\varphi)/A(\psi) = D > 1$, we know $A(\psi)$ is strictly less than $A(\varphi)$. Since $A(\varphi) < 4\pi$,  we can calculate that $g = 0$ or $1$ from $A(\psi) < 4\pi$.
	
	We know that $\orb$ cannot have $g=0$ and $m=3$ since triangular orbifolds are excluded in  flexible coverings. Since  we showed $3-2g \leq m \leq 5-4g$, the only remaining options for $g$ and $m$ are $g=0$ and $m = 4$, $g=0$ and $m=5$, or $g=1$ and $m=1$.
\end{proof}

Next we put  upper  and lower bounds on the degree of the flexible cover.

\begin{lemma}\label{Inequality on D} The degree $D$ of $p$ is at least 3 and most $12$.
\end{lemma}
\begin{proof}We first  establish $3\leq D$. Since $D \neq 1$, we must rule out $D=2$. If $D = 2$, then the local degree of each point on $(S,\varphi)$ is either 1 or 2. However, $(S,\varphi)$ must contain a cone point $y$. Since, $\Theta_y >2\pi$, the local degree of $y$ being at most 2  implies $\Theta_{p(y)} > \pi$, which is impossible because $\orb$ is an orbifold. Thus, $D \geq 3$.
	
	We now  show $D \leq 12$. 
	As in the proof of Lemma \ref{bounds on orbifold points}, $$D \leq \frac{4\pi}{-2\pi(2-2g-\sum_{i=1}^m(1-\frac{1}{r_i}))},$$ where $(g;r_1,\dots,r_m)$ is the signature of $\orb$.  We can minimize the denominator to get the corresponding upper bound for $D$ by determining the smallest possible value for $A(\psi)$. 
	
	Since $A(\psi)$ will increase with genus and the order of the orbifold points,  the smallest area orbifold will have genus $g =0$ and $m=4$ (we can rule out $m=3$ as $\orb$ must be non-triangular for a flexible cover).
	Note that $A(\psi) = 0$ corresponds to the signature $(0; 2,2,2,2)$ which is not a hyperbolic orbifold. Thus the smallest area \textit{hyperbolic} orbifold $\orb$ is $(0; 2,2,2, 3)$ with $A(\psi) = \pi/3$. Therefore,
	$D \leq \frac{4\pi}{\pi/3}  = 12. $
\end{proof}

We are now ready to compile our list of possible signatures. Recall that Lemma \ref{lem:order(x)<D} says that for a given degree $D$, each $r_i \leq D$. In particular, for a fixed $g$, $m$, and $D$, there are only finitely many signatures possible for $\orb$. The orbifold must also satisfy $A(\varphi)/A(\psi) = D$. As $A(\varphi) < 4\pi$, this means $D < 4\pi / A(\psi)$.
Hence we can directly compute all of the signatures and degrees $D$ so that  
\begin{enumerate}
	\item $D \in \{3,4,5 \dots, 12\}$ (Lemma \ref{Inequality on D})
	\item\label{item:g_and_m} $(g,m) \in \{(0,4),(0,5),(1,1)\}$ (Lemma \ref{bounds on orbifold points});
	\item\label{item:order<D} each $r_i \leq D$ (Lemma \ref{lem:order(x)<D});
	\item\label{item:odd} if $r_i$ is odd, then $r_i$ divides $D$ (Lemma \ref{lem:odd|D});
	\item\label{item:area} $\displaystyle  D < \frac{4\pi}{ 2\pi\left( 2g-2 + \sum_{i=1}^m \left(1-\frac{1}{r_i}\right) \right)} $.
\end{enumerate}
Table \ref{table:list} contains all of the signature and degree combination that satisfy the above requirements. This list was generated using computer code out of convenience, but could also be done by hand.

\begin{table}[h]
	\begin{tabular}{|c|c|c|}
		\hline
		\multicolumn{1}{|c|}{\textbf{Signature}} & \multicolumn{1}{c|}{\textbf{Degree}} & \textbf{Area} \\ \hline
		(0; 2, 2, 2, 3)    & 3, 6, 9 & $\pi/3$  \\ \hline
		(0; 2, 2, 3, 3)    & 3       & $2\pi/3$ \\ \hline
		(0; 2, 3, 3, 3)    & 3       & $\pi$    \\ \hline
		(0; 2, 2, 2, 4)    & 4, 5, 6, 7    & $\pi/2$  \\ \hline
		(0; 2, 2, 2, 5)    & 5       &  $3\pi/5$                     \\ \hline
		(0; 2, 2, 2, 2, 2) & 3       & $\pi     $\\ \hline
		(1; 2) &     3   &   $\pi$   \\ \hline
	\end{tabular}
	\caption{List of possible signatures and degrees for a flexible covering of a genus 2 surfaces to a hyperbolic orbifold.}\label{table:list}
\end{table}

We conclude by showing  some of the orbifolds from Table \ref{table:list} do not admit flexible coverings by any genus 2 negatively curved hyperbolic cone surface.

\begin{lemma}\label{lem:too_many_even_order}
	Suppose $p \colon (S,\varphi) \to \orb$ is a flexible covering map of  finite degree $D$. The following situations are impossible:
	\begin{itemize}
		\item $D=9$ and $\orb $ has signature $(0; 2, 2, 2, 3)$, 
		\item $D = 5$ or $7$ and $\orb $ has signature $(0; 2, 2, 2, 4)$,
		\item $D=5$ and $\orb $ has signature $(0; 2, 2, 2, 5)$, 
		\item $D=3$ and and $\orb $ has signature $(0; 2, 2, 2, 2, 2)$.
	\end{itemize}
\end{lemma} 

\begin{proof}
	For all of these cases, we  use   Lemma \ref{n_max} to establish an upper bound on the number of possible cone points on $(S,\varphi)$. We then see that for all of these orbifolds, the number of cone points on $(S,\varphi)$ is less than the number of even order orbifold points on $\orb$. This ensures that at least one even order orbifold point $x \in \orb$ has a preimage $p^{-1}(x)$ containing only regular points of $(S,\varphi)$. This implies that $D$ is even by  Lemma \ref{regular preimages}. However, in all of these cases, this creates a contradiction as $D$ is odd.
	
	To illustrate this line of reason, consider the example of signature $(0; 2, 2, 2, 3)$ and $D = 9$. 
	If  there exists a flexible covering map $p\colon (S,\varphi) \to \orb$,  then $(S,\varphi)$ has at most 1 cone point by Lemma \ref{n_max}. Since $\orb$ has 3 even order orbifold points, there exist an even order orbifold point $x \in \orb$ where every point of $p^{-1}(x)$  is a regular point of $(S,\varphi)$.  Lemma  \ref{regular preimages} then say  $D$ must be even, contradicting  $D = 9$.

    The remaining cases are similar.
\end{proof}

\section{Cone point and local degree restriction of flexible coverings}\label{sec:local_degree_restrictions}
We now examine the signatures in Table \ref{table:list} that are not prohibited from having a   flexible cover $p\colon (S,\varphi) \to \orb$  by Lemma \ref{lem:too_many_even_order}. We  determine what requirements the existence of a flexible cover  imposes on the  cone points of $(S,\varphi) $ and the local degrees of the preimages of orbifold points.  This will verify Theorem \ref{intro_thm:orbifolds_and_degrees} from the introduction.

Since the branch points of a flexible cover are exactly the orbifold points, these local degree restriction  also place restriction on the $\sym(D)$-tuple $\bt_p$ for the permutation representation of $p$.  In what follows, the statements about the tuple $\bt_p$ follows from the applying Theorem \ref{theorem:banched_cover_permutation} to the local degree restriction established in the lemmas.

\begin{lemma}\label{lem:degree_3_genus_0}
	Let $\varphi \in \hyp(S)$ and $\orb$ be an orbifold with signature $(0;2,2,2,3)$, $(0;,2,2,3,3)$, or $(0;2,3,3,3)$.  If there exists a degree 3 flexible cover $p \colon (S,\varphi) \to \orb$, then   $\pld(x) = (3)$ for each orbifold point $x \in \orb$. 
	Further, if $\order(x) = 2$, then $y=p^{-1}(x)$ is a cone point with $\Theta_y = 3\pi$.
    This implies that $\bt_p = [s_1,s_2,s_3,s_4]$ where $s_1s_2s_3 = s_4^{-1}$ and each $s_i$ is a 3-cycle in $\sym(3)$.
\end{lemma}

\begin{proof}
	Suppose there is a degree 3 flexible cover $p \colon (S,\varphi) \to \orb$. Let $x\in \orb$ be an orbifold point and  $y \in p^{-1}(x)$. Since $\Theta_y = \Theta_x \deg(y)$ must be at least $2\pi$, $\deg(y) \neq 1$.  Equation \ref{eq:sum_local_degree}  therefore forces $p^{-1}(x)$ to  contain only one point $y$ with $\deg(y) =3$.     
	If $\order(x) = 3$, then $y$ must be a regular point because $p$ is flexible.	If $\order(x) = 2$, then $y$ is a cone point with $\Theta_y = 3\pi$.
\end{proof}

\begin{lemma}\label{lem:degree_3_genus_1}
	Let $\varphi \in \hyp(S)$ and $\orb$ be an orbifold with signature $(1;2)$. Let $x \in \orb$ be the orbifold point. If there exists a degree $3$ flexible cover $p\colon (S,\varphi) \to \orb$, then  $\cone(\varphi) = (3\pi)$ and $\pld(x) = (3)$.   This implies $\bt_p = [s_1,s_2,s_3]$ where $s_1s_2s_1^{-1}s_2^{-1} = s_3$ and $s_3$ is a 3-cycle.    
\end{lemma}

\begin{proof}
	Suppose there is a degree 3 flexible cover $p \colon (S,\varphi) \to \orb$. Let $x\in \orb$ be the orbifold point. Since $p$ is flexible, $p^{-1}(x)$ must contain all the cone points of $(S,\varphi)$. Let $y \in p^{-1}(x)$ be a cone point. Now $\deg(y) >2$ as $\Theta_y > 2\pi$ and $\Theta_x = \Theta_y / \deg(y)$. Since $p$ has degree 3, the only option is $\deg(y) = 3$. This implies $y$ is the only point in $p^{-1}(x)$ and $\Theta_y = 3\pi$.  Hence $\cone(\varphi) = (3\pi)$ and $\pld(x) = (3)$.
\end{proof}

\begin{lemma}\label{lem:degree_4}
	Let $\varphi \in \hyp(S)$ and $\orb$ be an orbifold with signature $(0;2,2,2,4)$. If there exists a degree 4 flexible cover $p \colon (S,\varphi) \to \orb$, then $\cone(\varphi)= (4\pi)$ and 
	\begin{itemize}
		\item one of the orbifold points, $x_2^1$, of order 2, has $\pld(x_2^1) = (4)$ and $p^{-1}(x_2^1)$ is the single cone point of $(S,\varphi)$;
		\item two of the orbifold points, $x_2^2$, $x_2^3$,  of order $2$ have $\pld(x_2^i) = (2,2)$;
		\item the order 4 orbifold point $x_4$  has $\pld(x_4) = 4$.
	\end{itemize}
	This implies $\bt_p = [s_1,s_2,s_3,s_4]$ where $s_1s_2s_3=s_4^{-1}$, $s_4$ is a 4-cycle, and $[s_1,s_2,s_3]$ contains one element with cycle structure $(4)$ and two elements with cycle structure $(2,2)$.
\end{lemma}

\begin{proof}
	Let $y$ be a cone point on $\varphi$ and $x_4$ be the orbifold point with order 4. If $p(y) = x_4$,  then  $$ \deg(y)\frac{\pi}{2} = \deg(y) \Theta_{x_4} =\Theta_y > 2\pi \implies \deg(y) >4.$$ As $\deg(y)$ cannot be larger than the degree of $p$, this implies $p(y) \neq x_4$
	Therefore, $p(y)$ must be an order 2 orbifold point of $\orb$. Let $x_2^1 = p(y)$.
	
	Since $\Theta_y >2\pi$ and $\Theta_y = \deg(y) \pi$, $\deg(y)$ is either $3$ or $4$. If $\deg(y) =3$, then $p^{-1}(x_2^1)$ must contain a point $z$ with $\deg(z) = 1$. But this would imply $\Theta_z = \pi$, which would contradict $\varphi$ being negatively curved. Hence $\deg(y)$ must be $4$, $\Theta_y = 4\pi$, and $\pld(x_2^1) = (4)$.
	
	Now, the only way to satisfy the formula $ A(\varphi)/A(\psi) = 4$ when $(S,\varphi)$ contains at least 1 cone point of angle $4\pi$ is for $\cone(\varphi) = (4\pi)$.
	Since $(S,\varphi)$ contains only one cone point,  every other orbifold point of $\orb$ has only regular points in their preimage. This implies $\pld(x_2^2) = \pld(x^3_3) = (2,2)$ and $\pld(x_4) = (4)$ where $x_2^2$, $x_2^3$ are the other two  order $2$ orbifold points.
\end{proof}

\begin{lemma}\label{lem:degree_6_2,2,2,4}
	Let $\varphi \in \hyp(S)$ and $\orb$ be an orbifold with signature $(0;2,2,2,4)$. If there exists a degree 6 flexible cover $p \colon (S,\varphi) \to \orb$, then  
	\begin{itemize}
		\item $(S,\varphi)$ has a single cone point $y$ with $\Theta_y = 3\pi$;
		\item if $x \in \orb$ is the order $4$ orbifold point, then $p^{-1}(x) = y$ and $\pld(x) = (6)$;
		\item if $x\in \orb$ is an order 2 orbifold points then $\pld(x) = (2,2,2)$.
	\end{itemize}
	This implies $\bt_p = [s_1,s_2,s_3,s_4]$ where  $s_1s_2s_3=s_4^{-1}$, $s_4$ is a 6-cycle, and each of $s_1,s_2,s_3$ has cycle structure $(2,2,2)$.
\end{lemma}

\begin{proof}
	Let $x \in \orb$ be the order 4 orbifold point. We first show that $p^{-1}(x)$ cannot contain any regular points of $(S,\varphi)$. For each $y \in p^{-1}(x)$, $\deg(y) = \Theta_y/\Theta_{x}$. If $y$ is a regular point, then $\deg(y) = 4$ as $\Theta_y = 2\pi$ and $\Theta_{x_4} = \pi/2$. Since the local degrees for point in $p^{-1}(x)$ must sum to 6, this would imply that $p^{-1}(x)$ contains at least one more point $z$ with $\deg(z) \leq 2$. However, this would imply $$\frac{\pi}{2} \deg(z) = \Theta_{z} <\pi,$$ which contradicts $\varphi$ being negatively curved. Thus, $p^{-1}(x)$ can only contain cone points of $(S,\varphi)$.
	
	Next we show $p^{-1}(x)$ must be only a single cone point. Let $y \in p^{-1}(x)$. By the proceeding paragraph $\Theta_y >2\pi$. Since $\frac{\pi}{2} \deg(y) = \Theta_y$, we have $\deg(y) > 4$. Repeating the previous argument about points with local degree at most 2 shows that $\deg(y) \neq 5$. Hence $\deg(y) = 6$ and $p^{-1}(x) = y$.
	
	The above shows that $(S,\varphi)$ has at least one cone point with angle $\frac{\pi}{2} \cdot 6 = 3\pi$. Thus, the only way to satisfy the formula $A(\varphi)/A(\psi) = 6$ is for $(S,\varphi)$ to have exactly one cone point. Thus $\cone(\varphi) = (3\pi)$.
	
	Since $(S,\varphi)$ contains only one cone point, if $x \in \orb$ is an order 2 orbifold point, then $p^{-1}(x)$ must contain only regular points. Hence, $\pld(x) = (2,2,2)$.
\end{proof}

\begin{lemma}\label{lem:degree_6}
	Let $\varphi \in \hyp(S)$ and $\orb$ be an orbifold with signature $(0;2,2,2,3)$. If there exists a degree 6 flexible cover $p \colon (S,\varphi) \to \orb$, then one of the following holds.
	\begin{enumerate}
		\item $(S,\varphi)$ has two cone points, $y,y'$, each with cone angle $3\pi$ and with $p(y) = p(y')$. If $x_2^1,x_2^2,x_2^3 \in \orb$ are the order 2 orbifold points where $p(y) = x_2^1$, then $\pld(x_2^1)  = (3,3)$ and $\pld(x_2^2) = \pld(x_2^3) = (2,2,2)$. If $x_3 \in \orb$ is the order 3 orbifold point, then $\pld(x_3) = (3,3)$. This implies $\bt_p = [s_1,s_2,s_3,s_4]$ where $s_1s_2s_3 = s_4^{-1}$, $s_4$ has cycle structure $(3,3)$, and $[s_1,s_2,s_3]$ contains one element with cycle structure $(3,3)$ and two elements with cycle structure $(2,2,2)$.
		\item $(S,\varphi)$ has one cone point, $y$,  with cone angle $4\pi$. If $x_2^1,x_2^2,x_2^3 \in \orb$ are the order 2 orbifold points where $p(y) = x_2^1$, then $\pld(x_2^1) = (2,4)$, $\pld(x^2_2)=\pld(x^3_3) = (2,2,2).$  If $x_3 \in \orb$ is the order 3 orbifold point, then $\pld(x_3) = (3,3)$. This implies $\bt_p = [s_1,s_2,s_3,s_4]$ where $s_1s_2s_3 = s_4^{-1}$, $s_4$ has cycle structure $(3,3)$, and $[s_1,s_2,s_3]$ contains one element with cycle structure $(2,4)$ and two elements with cycle structure $(2,2,2)$.
	\end{enumerate}
	
\end{lemma}

\begin{proof}
	Since the degree of $p$ is 6, the local degree at each cone point of $(S,\varphi)$ is in $\{3,4,5,6\}$.
	If there is a cone point $y \in (S,\varphi)$ with $\deg(y) = 6$, then $\Theta_y = 6\pi$ because $\Theta_{p(y)} = \pi$. Since adding cone points to a negatively curved hyperbolic metric will lower the area, $\Theta_y = 6\pi$ implies $$ A(\varphi) \leq -2\pi \left( -2 - (1-3)\right) = 0.$$ However, this would contradict the fact that $\varphi$ is a hyperbolic metric, so $\deg(y) \neq 6$. If instead there is a cone point  $y$ with $\deg(y) =5$, then by Equation \ref{eq:sum_local_degree}, there is exactly one other point $z \in (S,\varphi)$ with $p(z) = p(y)$ and $\deg(z) = 1$. Since $\Theta_{p(y)} = \pi$, this would mean $\Theta_z = \pi$ as well, but this is a contradiction to the fact that $\varphi$ is negatively curved. Hence each cone point $y$ has $\deg(y) \in \{ 3,4\}$.
	
	By Lemma \ref{n_max}, the maximum number of cone points for $(S,\varphi)$ is $n_{max} = 2$. We now argue in cases depending on the number of cone points of $(S,\varphi)$.
	
	\textbf{Case 1:} Suppose $(S,\varphi)$ has two cone points. Using the formula $ A(\varphi)/A(\psi) = 6,$ we see that this case only occurs when each of the cone points has angle $3\pi$, thus $\cone(\varphi) = (3\pi, 3\pi)$. We now verify the local degrees at the preimages of the orbifold points.
	
	First consider the orbifold point $x_3$ of order $3$. Since $p$ is flexible, $p^{-1}(x_3)$ must contain only regular points of $(S,\varphi)$. This implies that the local degree at each point in $p^{-1}(x_3)$ is $3$, which means $\pld(x_3) = (3,3)$.
	
	Now consider consider the orbifold points, $x_2^1,x_2^2,x_2^3$, of order $2$. Since $p$ is flexible, at least one of these points has a cone point in its preimage. Without loss of generality, assume $p^{-1}(x_2^1)$ contains a cone point $(S,\varphi)$.  Let $\{y_1,\dots, y_k\} = p^{-1}(x_2^1)$ where at least $y_1$ is a cone point. Since $\Theta_{y_1} = 3\pi$, $\deg(y_1) = 3$. Thus $\sum_{j=2}^k \deg(y_j) = 3$. Since $\varphi$ is negatively curved, none of the $y_j$ can have local degree $1$. Moreover, if any $y_j$ is a regular point, then $\deg(y_j) = 2$. Thus, the only possibility is that $k=2$ and both $y_1$ and $y_2$ are cone points.  Thus $p(y_1) = p(y_2) = x_2^1$ and $\pld(x_2^1) = (3,3)$. Since $y_1$ and $y_2$ are the only cone points of $(S,\varphi)$, both $p^{-1}(x_2^2)$ and $p^{-1}(x_2^3)$ contain entirely regular points. Since each of these regular points must have a local degree of $2$,  $\pld(x_2^2) = \pld(x_2^3) = (2,2,2)$.

	\textbf{Case 2:} Suppose   $(S,\varphi)$ has one cone point $y$. Let $x_2^1 = p(y)$. Since $\deg(y) \in \{3,4\}$, the preimage $p^{-1}(x_2^1)$ must contain at least one additional points besides $y$. Since $(S,\varphi)$ has only one cone point, every point in $p^{-1}(x_2^1) - \{y\}$ is a regular point and each of these regular points has local degree $2$. Thus, the only possibility is that $p^{-1}(x_2^1)$ contains two points---the cone point $y$ and a regular point $z$---and that $\deg(y) =4$. Thus, $\cone(\varphi) = (4\pi)$ and $\pld(x_2^1) = (2,4)$. Just as in Case 1, we now know that the preimages of all other orbifold points each contain only regular points. Thus,we can conclude that $\pld(x_2^2)=\pld(x_2^3) = (2,2,2)$ and $\pld(x_3) = (3,3),$ where $x_i^j$ are the other orbifold points of order $i$.
\end{proof}

\section{Existence  of flexible covers}
The work in Section \ref{sec:local_degree_restrictions} gave restrictions on the behavior of the local degrees of flexible coverings by a genus 2 surface. Using Theorem \ref{theorem:banched_cover_permutation}, this also gives restrictions on the cycle structure of the permutation representation of the covering map. In this section, we will show that every permutation representation that meets these restrictions produces a flexible covering by a genus 2 surface.

\subsection{The $(1;2)$ orbifold}
Let $\orb$ be an orbifold with signature $(1;2)$ and let $x \in \orb$ be the orbifold point. As in Section \ref{subsec:permutations}, let $\overline{O}$ be the surface with one boundary components obtained by deleting a small open disk centered at the orbifold point $x$. Let $\delta$ be the boundary of the open disk around $x$, and let $\alpha,\beta, \gamma$ be the loops show in Figure \ref{fig:1-holed_torus}. Hence, $$ \pi_1(\overline{O}) = \langle \alpha,\beta,\gamma \mid \alpha \beta \alpha^{-1} \beta^{-1} = \gamma \rangle.$$ 

As described at the end of Section \ref{subsec:permutations}, there is a one-to-one correspondence between homomorphisms $\Phi \colon \pi_1(\overline{O}) \to \sym(3)$ and tuples, $[s_1,s_2,s_3]$, of elements of $\sym(3)$ that satisfy $s_1s_2 \overline{s}_1 \overline{s}_2 = s_3$. Specifically, such a tuple corresponds to the homomorphism determined by $\alpha \to s_1$, $\beta \to s_2$, and $\gamma \to s_3$.
By Theorem \ref{theorem:banched_cover_permutation}, each tuples also determines a degree $3$  branched covering map of some surface onto $\orb$ with conjugacy classes of these tuples correspond to equivalence classes of degree $3$ branched covers. We can therefore enumerate the flexible covers of a $\orb$ by finding all conjugacy classes of tuples that satisfy these requirements.

\begin{proposition}\label{prop:deg_3_genus_1}
	Let $\orb$ be an orbifold with signature $(1;2)$ and consider the tuples
	$$[(021), (01), (012)] \text{ and }   [(01),(021),(021)].$$
	Each of these  tuples corresponds to a topological branched covering map $p \colon S \to \orb$ so that if $\varphi$ is the pull back of the hyperbolic metric on $\orb$, then  $\varphi \in \hyp(S)$ and $p\colon (S,\varphi) \to \orb$ is a degree 3 flexible cover. Moreover, for all $\varphi \in \hyp(S)$, every degree 3 flexible covering of $(S,\varphi) \to \orb$ is covering space equivalent to exactly one of these two flexible coverings.
\end{proposition}

\begin{proof}
	Let $[s_1,s_2,s_3]$ be an tuple in $\sym(3)$ so that $s_1s_2\overline{s}_1\overline{s}_2 = s_3$, $s_3$ is a 3-cycle, and $[s_1,s_2]$ contains a 3-cycle and a 2-cycle. By Theorem \ref{theorem:banched_cover_permutation}, any such tuple produces a surface $Z$ and a branched cover $p \colon Z \to \orb$. 
	
	Let $\varphi$ be the hyperbolic metric on $Z$ obtained by pulling back the hyperbolic metric on $\orb$ under $p$. By Theorem \ref{theorem:banched_cover_permutation}, the preimage of the orbifold point $x \in \orb$ is a single point $y \in Z$ with $\Theta_y = \deg(y) \Theta_x = 3\pi$. Hence $\varphi \in \hyp(Z)$ and $p$ is a flexible covering map. We verify that $Z$ has genus $h=2$ using the fact that $A(\phi)/A(\psi) = 4$: 
	\begin{align*}
		\frac{-2\pi\left(2-2h + \frac{1}{2}\right)}{\pi} &= 3\\
		-2(5-4h) &= 6 \\
		h&=2.        
	\end{align*}

	Theorem \ref{theorem:banched_cover_permutation} and Lemma \ref{lem:degree_3_genus_1} say that for any $\varphi \in \hyp(S)$, any flexible cover $p \colon (S,\varphi) \to \orb$ must correspond to a tuple of elements in $\sym(3)$ with the conditions outlined in the first paragraph. Any tuple $[s_1,s_2,s_3]$ where $s_1s_2\overline{s}_1\overline{s}_2 = s_3$ and $s_1$ is a 3-cycle is conjugate to one  of $$ [(021), (01), (012)], \ [(021), (02), (012)], \ \text{or} \ [(021), (12), (012)].$$ By conjugating by either $(021)$ or $(012)$ we see that all three of these tuple are in fact conjugate. Similar reasoning applies to show that any tuple where $s_2$ is a 3-cycle  is conjugate to $[(01),(021),(021)]$.
\end{proof}

\subsection{Genus 0 orbifolds} Let $\orb$ be an orbifold with signature $(0;r_1,r_2,r_3,r_4)$ and let $x_i \in \orb$ be the orbifold point with order $r_i$. As in Section \ref{subsec:permutations}, let $\overline{O}$ be the surface with four boundary components obtained by deleting a small open disk centered at each orbifold point $x_i$. Let $\delta_i$ be the boundary of the open disk around $x_i$, and let $\gamma_i$ be the loops show in Figure \ref{fig:4-holed_sphere}. Thus, we have $$ \pi_1(\overline{O}) = \langle \gamma_1,\gamma_2,\gamma_3,\gamma_4 \mid \gamma_1\gamma_2\gamma_3 =\gamma_4^{-1} \rangle.$$ 

As described at the end of Section \ref{subsec:permutations}, there is a one-to-one correspondence between homomorphisms $\Phi \colon \pi_1(\overline{O}) \to \sym(D)$ and 4-tuples of elements, $[s_1,s_2,s_3,s_4]$, of $\sym(D)$ that satisfy $s_1s_2s_3 = s_4^{-1}$. Namely, the tuple $[s_1,s_2,s_3,s_4]$ corresponds with the homomorphism defined by sending each $\gamma_i$ to $s_i$.
By Theorem \ref{theorem:banched_cover_permutation}, each tuple also determines a degree $D$ topological branched covering map of some surface onto $\orb$ with conjugacy classes of these tuples correspond to equivalence class of degree $D$ branched covers.

In Section \ref{sec:local_degree_restrictions}, we determined restrictions on local degrees of the preimages of orbifolds under any possible flexible covering by $S$.  Theorem \ref{theorem:banched_cover_permutation} says these restrictions on local degrees  translate into restrictions on the cycle structures of the elements $s_i$ in a tuple defining a homomorphism into $\sym(D)$. We  verify that if a tuple $[s_1,s_2,s_3,s_4]$ satisfies these restriction, then the corresponding covering map will produce a flexible metric on $S$ by pulling back the metric on $\orb$.

\begin{lemma}\label{lem:existance_deg_3_genus_0}
	Let $\orb$ be an orbifold with signature $(0;2,2,2,3)$, $(0;2,2,3,3,)$, or $(0;,2,3,3,3)$. Let $[s_1,s_2,s_3,s_4]$ be a tuple of elements of $\sym(3)$ so that $s_1s_2s_3 =s_4^{-1}$,  and each $s_i$ is a 3-cycle.  Any such tuples corresponds to a topological branched covering map $p \colon S \to \orb$ so that if $\varphi$ is the pull back of the hyperbolic metric on $\orb$, then  $\varphi \in \hyp(S)$ and $p\colon (S,\varphi) \to \orb$ is a degree 3 flexible cover.
\end{lemma}

\begin{proof}
	Let $[s_1,s_2,s_3,s_4]$ be a tuple of elements of $\sym(3)$ so that $s_1s_2s_3 =s_4^{-1}$,  and each $s_i$ is a 3-cycle. For any such tuple, Theorem \ref{theorem:banched_cover_permutation} provides a surface $Z$ and a topological branched cover $p \colon Z \to \orb$. 
	
	Let $\varphi$ be the metric on $Z$ obtained by pulling back the hyperbolic metric on $\orb$ under $p$. Theorem \ref{theorem:banched_cover_permutation} say that for each orbifold point $x_i \in \orb$, $p^{-1}(x_i)$ is a single point $y_i$ and $\deg(y_i) = 3$. Thus, $\Theta_{y_i} = 3\pi$ if $\order(x_i) = 2$ and $\Theta_{y_i}= 2\pi$ if $\order(x_i) = 3$. Since at least one $x_i$ has order $2$, $\varphi \in \hyp(S)$ and  $p \colon (Z,\varphi) \to \orb$ is a flexible covering.
	
	It remains to check that $Z$ is the genus 2 surface $S$.  We use the formulas $D=A(\varphi)/A(\psi)$ and $\Theta_{y_i} = 6\pi/\order(x_i)$ to solve for the genus, $h$ of $Z$. 
	\begin{align*}
		\frac{2-2h-\sum_{i=1}^4\left(1-\frac{3}{\order(x_i)}\right)}{ 2 - \sum_{i=1}^4  \left(1-\frac{1}{\order(x_i)}\right)} &= 3   \\
		2-2h - \sum_{i=1}^4\left(1-\frac{3}{\order(x_i)}\right) &= 6 - \sum_{i=1}^4  \left(3-\frac{3}{\order(x_i)}\right) \\
		-2h -2 - \sum_{i=1}^4\frac{3}{\order(x_i)} &= -6 - \sum_{i=1}^4  \frac{3}{\order(x_i)} \\
		h &= 2
	\end{align*}
\end{proof}

\begin{lemma}\label{lem:existance_deg_4}
	Let $\orb$ be a  $(0;2,2,2,4)$ orbifold. Let $[s_1,s_2,s_3,s_4]$ be a tuple of elements of $\sym(4)$ so that  \begin{itemize}
		\item $s_1s_2s_3 =s_4^{-1}$,
		\item  $s_4$ has cycle structure $(4)$,
		\item    one of $s_1,s_2,s_3$ has cycle structure $(4)$ and the other two  have cycle structure $(2,2)$.
	\end{itemize} Every such tuple corresponds to a topological branched covering map $p \colon S \to \orb$ so that if $\varphi$ is the pull back of the hyperbolic metric on $\orb$, then  $\varphi \in \hyp(S)$ and $p\colon (S,\varphi) \to \orb$ is a degree 4 flexible cover.
\end{lemma}

\begin{proof}
	Let $[s_1,s_2,s_3,s_4]$ be any tuple of elements of $\sym(4)$ that is describe in Lemma \ref{lem:existance_deg_4}.
	For any such tuple, Theorem \ref{theorem:banched_cover_permutation} provides a surface $Z$ and a topological branched cover $p \colon Z \to \orb$. 
	
	Let $t_1$ be the $4$-cycle  term of $[s_1,s_2,s_3]$  and $t_2,t_3$ be the two terms with cycle structure $(2,2)$. Let $t_4=s_4$. There is now a permutation $\sigma \in \sym(4)$ so that $t_i = s_{\sigma(i)}$. For each $i \in \{1,2, 3,4\}$, define $o_i = x_{\sigma(i)}$.
	
	Let $\varphi$ be the metric on $Z$ obtained by pulling back the hyperbolic metric on $\orb$ under $p$.
	By Theorem \ref{theorem:banched_cover_permutation}, $p^{-1}(o_4)$ is a single point $y_4$ with $\Theta_{y_4} = 2\pi$, $p^{-1}(o_1)$ is a single point $y_1$ with $\Theta_{y_1} = \deg(y_1) \pi = 4\pi$, and $p^{-1}(o_2)$ and $p^{-1}(o_3)$ each contain two regular points of $(Z,\varphi)$. Hence $\varphi$ is a hyperbolic cone metric on $Z$ and $p \colon (Z,\varphi) \to \orb$ is a flexible covering.
	
	To check that $Z$ is the genus 2 surface $S$, we use the formula $4=A(\varphi)/A(\psi)$ to solve for the genus of $Z$, and see that it is 2.     \end{proof}

\begin{lemma}\label{lem:existance_deg_6_2,2,2,4}
	Let $\orb$ be a $(0;2,2,2,4)$ orbifold.
	Let $[s_1,s_2,s_3s_4]$ be a tuple of  of elements of $\sym(6)$ so that $s_1s_2s_3 =s_4^{-1}$,  $s_1,s_2,s_3$ all have cycle structure $(2,2,2)$, and $s_4$ has cycle structure $(6)$. Any such tuple corresponds to a topological branched covering map $p \colon S \to \orb$ so that if $\varphi$ is the pull back of the hyperbolic metric on $\orb$, then  $\varphi \in \hyp(S)$ and $p\colon (S,\varphi) \to \orb$ is a degree 6 flexible cover.
\end{lemma}

\begin{proof}
	
	Let $[s_1,s_2,s_3,s_4]$ be any tuple of elements of $\sym(6)$ so that $s_1s_2s_3 =s_4^{-1}$,  $s_1,s_2,s_3$ have cycle structure $(2,2,2)$, and $s_4$ has cycle structure $(6)$. For any such tuple, Theorem \ref{theorem:banched_cover_permutation} provides a surface $Z$ and a topological branched cover $p \colon Z \to \orb$. 
	
	Let $\varphi$ be the metric on $Z$ obtained by pulling back under $p$ the hyperbolic metric on $\orb$.  For $i \in \{1,2,3\}$, Theorem \ref{theorem:banched_cover_permutation} says  $p^{-1}(x_i) = \{y_i^1,y_i^2,y_i^3\}$ and $\deg(y_i^j) = 2$ for each $i,j \in \{1,2,3\}$.  Since $\Theta_{y_i^j} = \deg(y_i^j) \Theta_{x_i} = 2\pi$, each $y_i^j$ is a regular point. Similarly, $p^{-1}(x_4)$ is a single point $y$ with $\deg(y) = 6$. Thus $\Theta_y = 3\pi = \Theta_{x_4} \deg(y)$. Hence $\varphi$ is a hyperbolic cone metric on $Z$ and $p \colon (Z,\varphi) \to \orb$ is a flexible covering.
	
	To check that $Z$ is the genus 2 surface $S$, we solve the formula $6=A(\varphi)/A(\psi)$ to solve for the genus of $Z$ and find that it is 2. 
\end{proof}

\begin{lemma}\label{lem:existance_deg_6_2,2,2,3}
	Let $\orb$ be an $(0;2,2,2,3)$ orbifold, and consider a tuple $[s_1,s_2,s_3,s_4]$ of elements of $\sym(6)$ where $s_1s_2s_3= s_4^{-1}$ and $s_4$ has cycle structure $(3,3)$. If two of $s_1,s_2,s_3$ have cycle structure $(2,2,2)$ and the third has cycle structure  $(2,4)$ or $(3,3)$, then $[s_1,s_2,s_3,s_4]$ corresponds to a topological branched covering map $p \colon S \to \orb$ so that if $\varphi$ is the pull back of the hyperbolic metric on $\orb$, then  $\varphi \in \hyp(S)$ and $p\colon (S,\varphi) \to \orb$ is a degree 6 flexible cover. 
\end{lemma}

\begin{proof}
	Let $[s_1,s_2,s_3,s_4]$ be a tuple of  elements of $\sym(6)$ where $s_1s_2s_3= s_4^{-1}$ and $s_4$ has cycle structure $(3,3)$. Assume two terms, $t_1$ and $t_2$,  of $[s_1,s_2,s_3]$ have cycle structure $(2,2,2)$. Let $t_3$ be the remaining term of $[s_1,s_2,s_3]$, and  let $t_4 = s_4$. There is $\sigma \in \sym(4)$ so that for each $i \in \{1,2,3,4\}$, $t_i = s_{\sigma(i)}$. Define $o_i = x_{\sigma(i)}$.
	
	Since $t_1s_4$ is an element of order $6$, the group generated by the $s_i$ is transitive. Thus, Theorem \ref{theorem:banched_cover_permutation} produces a connected surface $Z$ and topological branched cover $p \colon Z \to \orb$.  Let $h$ be the genus of $Z$ and $\varphi$ be the hyperbolic metric on $Z$ obtained by pulling back the hyperbolic metric on $\orb$.  For $i =1$ or $2$, Theorem \ref{theorem:banched_cover_permutation} plus Equation \ref{eq:cover_degree} says $p^{-1}(o_i)$ must contain three points $y_i^1,y_i^2,y_i^3$ each with $\Theta_{y_i^j}= 2\pi $. Similarly, $p^{-1}(o_4)$ contains two points $y_4^1,y_4^2$ each with $\Theta_{y_i^j} = 3 \cdot 2\pi/3 = 2\pi$.

	Suppose first that $t_3$ has cycle structure $(2,4)$.  Thus, $p^{-1}(o_3)$ contains two points  $y_3^1$ and $y_3^2$, where $\deg(y_3^1) =2$ and $\deg(y_3^2) = 4$. Thus, $\Theta_{y_3^1} = 2\pi$ and $\Theta_{y_3^2} = 4\pi$. Hence, $\varphi \in \hyp(Z)$ and $p \colon (Z,\varphi) \to \orb$ is a flexible covering.

	Now suppose that $t_3$ has cycle structure $(3,3)$.  In this case, $p^{-1}(o_3)$ contains two points each with cone angle $3\pi$. Thus $\varphi \in \hyp(Z)$ and  $\varphi \in \hyp(Z)$ and $p \colon (Z,\varphi) \to \orb$ is a flexible covering.    
	
	In both cases for $t_3$, the $\sum 1-\frac{\Theta_i}{2\pi}$ term of $A(\varphi)$ is equal to $1$. Thus, we can verify that $Z$ has genus $h=2$ with the calculation:
	\begin{align*}
		6 &= \frac{-2\pi(2-2h +1)}{\pi/3}\\
		6 & = -6(3-2h)\\
		h &= 2.
	\end{align*}
\end{proof}

\section{Counting Flexibility Classes}

We now count the $\MCG(S)$-orbits of the flexibility classes in $\hyp(S)$. According to Corollary \ref{cor:flex_class_iff_MCG-equiv}, we can count these flexibility classes by counting the signature equivalence classes of holonomy covers by $S$. We first count the number of signature equivalence classes of all flexible coverings by a genus two surface.

\begin{theorem}\label{thm:MCG_equivelence_of_flex_covers}
	Let $\orb$ be a hyperbolic orbifold. Suppose there is a metric $\varphi \in \hyp(S)$ and a degree $D$ flexible cover $p\colon (S,\varphi) \to \orb$. Let $\bt$ be the $\sym(D)$-tuple for $p$.  Exactly one of the following holds:
	\begin{enumerate}
		\item $D = 3$, the signature of $\orb$ is $(1;2)$, and $\bt$ is  signature equivalent to exactly one of $$ [(021),(01), (012) ] \text{ or } [(01),(12), (021)].$$
		\item $D=3$, the signature of $\orb$ is $(0;2,2,2,3)$ or $(0;2,3,3,3)$, and $\bt$ signature equivalent to $$ [(012),(012), (021), (021) ].$$
		\item $D=3$, the signature of $\orb$ is $(0;2,2,3,3)$, and $\bt$ is signature equivalent to one of $$ [(012),(012), (021), (021) ] \text{ or }  [(012), (021), (012), (021) ].  $$
		\item $D=4$,  the signature of $\orb$ is $(0;2,2,2,4)$,  $\bt$ signature equivalent to exactly one of 
		$$ [(0 1 2 3), (0 1)(2 3), (0 1)(2 3), (0 3 2 1)] \text{ or }
		[(0 1 2 3), (0 2)(1 3), (0 2)(1 3), (0 3 2 1)].$$
		
		\item $D=6$, the the signature of $\orb$ is $(0;2,2,2,4)$, and $\bt$ signature equivalent to $$[ (0 2)(1 5)(3 4), (0 1)(2 4)(3 5), (0 3)(1 5)(2 4), (0 2 5 3 4 1)].$$
		
		\item $D=6$, the the signature of $\orb$ is $(0;2,2,2,3)$, and $\bt$ signature equivalent to exactly one of 
		$$[(0 1 4 5)(2 3), (0 3)(1 5)(2 4), (0 4)(1 5)(2 3), (0 5 3)(1 2 4)] \quad \quad [(0 3 1)(2 4 5), (0 2)(1 5)(3 4), (0 2)(1 5)(3 4), (0 1 3)(2 5 4)]$$
		$$[(0 3 1)(2 4 5), (0 1)(2 4)(3 5), (0 1)(2 4)(3 5), (0 1 3)(2 5 4)] \quad \quad[(0 3 1)(2 4 5), (0 4)(1 5)(2 3), (0 4)(1 5)(2 3), (0 1 3)(2 5 4)].$$

	\end{enumerate}

\end{theorem}

\begin{proof}
	The limitation on the degree and signatures comes from the work in Section \ref{sec:local_degree_restrictions}. Section \ref{sec:local_degree_restrictions} also puts restriction on the cycle structure of the elements of $\bt$. We will show that in each case $\bt$ is signature equivalent to one of the listed tuples. The proof in each case has the same form. First we use the $\MCG\orb$ generators $H_{i,j}$, to assume that $\bt = [s_1,s_2,s_3,s_4]$ where each $s_i$ has a determined cycle structure. Then, we use computer code to generate a list of all conjugacy classes of $\sym(D)$-tuples that have that cycle structure and satisfy $s_1s_2s_3=s_4^{-1}$. This will give us a finite list of tuple where $\bt$ must be signature equivalent to one of the tuple on the list. To ensure that $\bt$ is signature equivalent to exactly one tuple on this list, we identify the $\MCG\orb$-orbits of the conjugacy classes of the $\sym(D)$-tuples.

	\textbf{Case: $D=3$ and signature $(1;2)$.} By Proposition \ref{prop:deg_3_genus_1}, $\bt$ must be conjugate to either $\bt_1=[(021), (01), (012)]$ or  $\bt_2=[(01),(021),(021)]$.  If $f$ is the order 2 homeomorphism of $\orb$ shown in Figure \ref{fig:order_2}, then $f \cdot \bt_1 = \bt_2$.  Hence $\bt$ is signature equivalent to $\bt_1$.

	\begin{figure}[ht]
		

             \def\svgwdith{2in}
             \hspace{1cm} 
\begingroup%
  \makeatletter%
  \providecommand\color[2][]{%
    \errmessage{(Inkscape) Color is used for the text in Inkscape, but the package 'color.sty' is not loaded}%
    \renewcommand\color[2][]{}%
  }%
  \providecommand\transparent[1]{%
    \errmessage{(Inkscape) Transparency is used (non-zero) for the text in Inkscape, but the package 'transparent.sty' is not loaded}%
    \renewcommand\transparent[1]{}%
  }%
  \providecommand\rotatebox[2]{#2}%
  \newcommand*\fsize{\dimexpr\f@size pt\relax}%
  \newcommand*\lineheight[1]{\fontsize{\fsize}{#1\fsize}\selectfont}%
  \ifx\svgwidth\undefined%
    \setlength{\unitlength}{63.54836495bp}%
    \ifx\svgscale\undefined%
      \relax%
    \else%
      \setlength{\unitlength}{\unitlength * \real{\svgscale}}%
    \fi%
  \else%
    \setlength{\unitlength}{\svgwidth}%
  \fi%
  \global\let\svgwidth\undefined%
  \global\let\svgscale\undefined%
  \makeatother%
  \begin{picture}(1,0.82783435)%
    \lineheight{1}%
    \setlength\tabcolsep{0pt}%
    \put(0.26366497,0.6175065){\color[rgb]{0,0,1}\makebox(0,0)[lt]{\lineheight{1.25}\smash{\begin{tabular}[t]{l}$\beta$\end{tabular}}}}%
    \put(-0.01567456,0.29958006){\color[rgb]{1,0,0}\makebox(0,0)[lt]{\lineheight{1.25}\smash{\begin{tabular}[t]{l}$\alpha$\end{tabular}}}}%
    \put(0,0){\includegraphics[width=\unitlength,page=1]{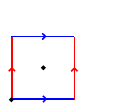}}%
    \put(0.70891412,0.68436599){\makebox(0,0)[lt]{\lineheight{1.25}\smash{\begin{tabular}[t]{l}$f$\end{tabular}}}}%
    \put(0,0){\includegraphics[width=\unitlength,page=2]{1_holed_torus_order_2_svg-tex.pdf}}%
  \end{picture}%
\endgroup%

		\caption{ The center point of the square is the orbifold point of the $(1;2)$ orbifold $\orb$. Flipping the square along the dashed diagonal is an order 2 homeomorphism that  exchanges $\alpha$ and $\beta$.}
		\label{fig:order_2}
	\end{figure}

	\textbf{Case: $D=3$ and  signature $(0;2,2,2,3)$, $(0;2,2,3,3,)$, or $(0;2,3,3,3)$.} By Lemma \ref{lem:degree_3_genus_0}, $\bt$ must be conjugate to a $\sym(3)$-tuple $[s_1,s_2,s_3,s_4]$ where $s_1s_2s_3=s_4^{-1}$ and each $s_i$ is a 3-cycle. There are three conjugacy classes of such tuples represented by 
	$$ [(0 1 2), (0 1 2), (0 2 1), (0 2 1)], \ 
	[(0 1 2), (0 2 1), (0 1 2), (0 2 1)], \text{ or }
	[(0 1 2), (0 2 1), (0 2 1), (0 1 2)].$$ 
	
	Since each of these tuples will generate cyclic group of order $3$,  the $\MCG\orb$-action is just permutation of the terms of the tuples. In the $(0;2,2,2,3)$ and $(0;2,3,3,3)$ case, we can use either $H_{2,1}$ and $H_{3,2}$ or $H_{3,2}$ and $H_{4,3}$ to establish that all three of these conjugacy classes are in the same $\MCG\orb$-orbit.
	Thus, $\bt$ will be signature equivalent to $[(0 1 2), (0 1 2), (0 2 1), (0 2 1)]$.
	
	In the $(0;2,2,3,3)$ case, the $\MCG\orb$-action can only exchange $s_1$ with $s_2$ and $s_3$ with $s_4$ because $\MCG\orb$ preserves the order of the cone points. Thus the conjugacy classes of $$ [(012),(012), (021), (021) ] \text{ and }  [(012), (021), (012), (021) ]  $$ are in separate $\MCG\orb$-orbits. Hence $\bt$ will be signature equivalent to exactly one of them.
	
	\textbf{Case: $D=4$ and signature $(0;2,2,2,4)$}. Using Lemma \ref{lem:degree_4} plus $H_{2,1}$ and $H_{3,1}$, we have that $\bt$ is signature equivalent to a tuple $[s_1,s_2,s_3,s_4]$ where $s_1$  is 4-cycle and $s_2,s_3$ are products of two 2-cycles, and $s_1s_2s_3 = s_4^{-1}$. Using the Python code in the appendix, we find three conjugacy classes of such tuples each represented by one of $$ \mathbf{t}_1= [(0 1 2 3), (0 1)(2 3), (0 1)(2 3), (0 3 2 1)],$$ $$ \mathbf{t}_2=  [(0 1 2 3), (0 1)(2 3), (0 3)(1 2), (0 1 2 3)],$$
	$$\mathbf{t}_3= [(0 1 2 3), (0 2)(1 3), (0 2)(1 3), (0 3 2 1)].$$ 
	
	Now $F_{2,1} \cdot \mathbf{t}_1  =  [(0321),(03)(12), (01)(23), (0321) ]$ which is conjugate to $\bt_2$ by the element $(13)$. 
	
	To see that $ \mathbf{t}_1$ is not signature equivalent to $\mathbf{t}_3$, observe that the group generated by the element of $\bt_3$ is the cyclic group generated by $(0123)$, while the group generated by $\bt_1$ is not cyclic. Since isomorphism class of the subgroup generated by the tuple is an invariant of signature equivalence class, $\bt_1$ and $\bt_3$ cannot be signature equivalent.
	
	\textbf{Case: $D=6$ and signature $(0;2,2,2,4)$.} By Lemma \ref{lem:degree_6_2,2,2,4},  $\bt$ is a tuple $[s_1,s_2,s_3,s_4]$ so that $s_1,s_2,s_3$ all have cycle structure $(2,2)$ and $s_4$ is a 4-cycle. Using the Python code in the appendix, there are exactly three conjugacy classes of such tuples represented by $$\bt_1 = [ (0 2)(1 5)(3 4), (0 1)(2 4)(3 5), (0 3)(1 5)(2 4), (0 2 5 3 4 1)] $$ 
	$$\bt_2 = [(0 2)(1 5)(3 4), (0 2)(1 3)(4 5), (0 5)(1 3)(2 4), (0 5 1 2 4 3)] $$ 
	$$\bt_3 = [ (0 2)(1 5)(3 4), (0 2)(1 3)(4 5), (0 3)(1 5)(2 4), (0 3 1 2 4 5)]. $$ 
	
	Now $H_{2,3} \cdot \bt_1 = [ (0 2)(1 5)(3 4), (0 3)(1 5)(2 4), (0 1)(2 4)(3 5), (0 2 5 3 4 1) ]$ which is conjugate to $\bt_2$ by the element $(041235)$. For $\bt_3$, we have $H_{1,3} \cdot \bt_1 = [ (0 3)(1 5)(2 4), (0 3)(1 4)(2 5), (0 2)(1 5)(3 4), (0 2 5 3 4 1)] $ which is conjugate to $\bt_3$ by the element $(15)(23)$.
	
	Since $\bt$ must be conjugate to one of $\bt_1,\bt_2,\bt_3$, $\bt$ is signature equivalent to $\bt_1$.
	
	\textbf{Case: $D=6$ and signature $(0;2,2,2,3)$.}
	By Lemma \ref{lem:degree_6}, $\bt$ is a tuple $[s_1,s_2,s_3,s_4]$ where $s_1s_2s_3=s_4^{-1}$, $s_4$ has cycle structure $(3,3)$, two of $s_1,s_2,s_3$ has cycle structure $(2,2,2)$, and the third has cycle structure either $(2,4)$ or $(3,3)$. Using $H_{2,1}$ and $H_{3,2}$, we can assume $s_1$ has cycle structure either $(2,4)$ or $(3,3)$ and $s_2,s_3$ have cycle structure $(2,2,2)$. 
	
	We will start be tackling the case where $s_1$ has cycle structure $(2,4)$. Using the Python code in the appendix, we find that there are two conjugacy classes the $\bt$ could belong to. These are represented by
	$$\bt_1 = [(0 1 4 5)(2 3), (0 3)(1 5)(2 4), (0 4)(1 5)(2 3), (0 5 3)(1 2 4)]$$
	and 
	$$\bt_2 = [ (0 1 4 5)(2 3), (0 4)(1 5)(2 3), (0 3)(1 5)(2 4), (0 5 3)(1 2 4)].$$
	
	\noindent Since $H_{3,2} \cdot \bt_1 = \bt_2$,  $\bt$ is signature equivalent to $\bt_1$.
	
	We now handle the $(3,3)$ case. Using the same Python code as before, there are nine conjugacy classes for $\bt$ in this case. Here are their representatives: 
	
	$$\bt_1 = [(0 3 1)(2 4 5), (0 2)(1 5)(3 4), (0 2)(1 5)(3 4), (0 1 3)(2 5 4)]$$ 
	$$\bt_2 = [(0 3 1)(2 4 5), (0 1)(2 4)(3 5), (0 1)(2 4)(3 5), (0 1 3)(2 5 4)]$$ 
	$$\bt_3 = [(0 3 1)(2 4 5), (0 1)(2 4)(3 5), (0 3)(1 2)(4 5), (0 3 2)(1 4 5)]$$
	$$\bt_4 = [(0 3 1)(2 4 5), (0 1)(2 4)(3 5), (0 2)(1 4)(3 5), (0 2 5)(1 3 4)]$$
	$$\bt_5 = [(0 3 1)(2 4 5), (0 1)(2 4)(3 5), (0 5)(1 2)(3 4), (0 5 1)(2 4 3)]$$
	$$\bt_6 = [(0 3 1)(2 4 5), (0 4)(1 5)(2 3), (0 1)(2 4)(3 5), (0 3 2)(1 4 5)]$$
	$$\bt_7 = [(0 3 1)(2 4 5), (0 4)(1 5)(2 3), (0 3)(1 5)(2 4), (0 1 4)(2 5 3)]$$
	$$\bt_8 = [(0 3 1)(2 4 5), (0 4)(1 5)(2 3), (0 4)(1 5)(2 3), (0 1 3)(2 5 4)]$$ 
	$$\bt_9 = [(0 3 1)(2 4 5), (0 4)(1 5)(2 3), (0 2)(1 4)(3 5), (0 3 1)(2 4 5)]$$
	
	We will verify that there are a total of three distinct signature equivalence classes among these nine tuple. First we show that each of the above tuples is signature equivalent to one of $\bt_1$, $\bt_2$, or $\bt_8$.

	\begin{itemize}
		\item  $F_{1,2}\cdot \bt_2 = [(0 1 5)(2 3 4), (0 5)(1 2)(3 4), (0 1)(2 4)(3 5), (0 1 3)(2 5 4) ] $, which is conjugate to $\bt_3$ by the element $(1325)$.
		\item $F_{1,3} \cdot \bt_2 = \bt_6$.
		\item $F_{2,4} \cdot \bt_2 = [(0 1 5)(2 3 4), (0 1)(2 4)(3 5), (0 4)(1 3)(2 5), (0 1 3)(2 5 4)]$ which is conjugate to $\bt_5$ by the element $(02)(14)$.
		\item $F_{1,2} \cdot \bt_4 = [ (0 1 5)(2 3 4), (0 5)(1 2)(3 4), (0 2)(1 4)(3 5), (0 2 5)(1 3 4)]$ which is conjugate to $\bt_5$ by the element $(1325)$.
		\item $H_{2,3} \cdot \bt_4 = [ (0 3 1)(2 4 5), (0 2)(1 4)(3 5), (0 1)(2 4)(3 5), (0 2 5)(1 3 4)]$ which is conjugate to $\bt_7$ be the element $(031)$.
		\item $F_{1,2}\cdot \bt_8 = [ (0 1 3)(2 5 4), (0 2)(1 4)(3 5) ,(0 2)(1 4)(3 5), (0 3 1)(2 4 5)]$ which is conjugate to $\bt_9$ by the element $(13)(24)$.
	\end{itemize}
	
	To verify that no two of $\bt_1$, $\bt_2$, and  $\bt_8$ are signature equivalent, consider the subgroup of $\sym(6)$ generated by these tuples. The isomorphism class of this subgroup is an invariant of the signature equivalence class of the tuple. Now, $\bt_1$ generates a cyclic group of order 6, $\bt_2$ generates a group of order 24, and $\bt_8$ generates a dihedral group of order 6. Since none of these are isomorphic, we have found three distinct signature equivalence classes of tuples. Hence $\bt$ is signature equivalent to exactly one  of them.
\end{proof}

Theorem \ref{thm:MCG_equivelence_of_flex_covers} say there are only 12 signature equivalence class of flexible covering in the genus 2 case. Thus there are at most 12 $\MCG(S)$-orbits of flexibility classes in $\hyp(S)$. However, it is possible that some of these flexible coverings  do not represent holonomy covers, making 12 an over count of the flexibility classes. Corollary \ref{cor:maximal_flexible_cover} showed that this over counting only occurs when the flexible cover is covered by a holonomy cover. The next lemma limits the situations where this nesting of flexible covers can occur.

\begin{lemma}\label{lem:non-degree_1_covers}
	Let $\orbone$ and $\orbtwo$ be  orbifolds so that for each $i \in\{1,2\}$ there is a degree $D_i$ flexible covering $p_i \colon (S,\varphi) \to \orbi$. If there exist a locally isometric branched cover $q \colon \orbone \to \orbtwo$ so that $p_2 = q \circ p_1$ and $q$ has degree at least $2$, then  $q$ has degree $2$ and  the signatures of $\orbone$ are $\orbtwo$ are respectively either
	\[(0;2,2,3,3) \text{ and } (0;2,2,2,3) \]
	or
	\[ (1;2) \text { and } (0;2,2,2,4).\]
\end{lemma}

\begin{proof}
	If there exists a degree $B$  locally isometric branched cover $q \colon \orbone \to \orbtwo$ so that $p_2 = q \circ p_1$, then $D_2 = B \cdot D_1$. According to Table \ref{table:list} and Lemma \ref{lem:too_many_even_order}, $D_1,D_2 \in \{3,4,6\}$. Thus, the only option is  $D_2 = 6$, $D_1 = 3$, and $B =2$.  Moreover, $\orbtwo$ must have signature $(0;2,2,2,4)$ or $(0;2,2,2,3)$.
	
	If $\orbtwo$ has signature $(0;2,2,2,4)$, then $A(\psi_1) = 2 A(\psi_2) = \pi$. Hence $\orbone$ has signature $(1;2)$ or $(0;2,3,3,3)$. Since none of the orbifold points of $\orbtwo$ have order divisible by $3$, it is impossible for $\orbtwo$ to be covered by an orbifold with signature $(0;2,3,3,3)$.
	
	If $\orbtwo$ has signature $(0;2,2,2,3)$, then $\orbone$ has signature $(0;2,2,3,3)$  as $A(\psi_1) = 2 A(\psi_2) = 2\pi/3$. 
\end{proof}

Using Lemma \ref{lem:non-degree_1_covers} and Corollaries \ref{cor:maximal_flexible_cover} and \ref{cor:flex_class_iff_MCG-equiv}, we can now count the signature equivalence classes of holonomy covers and by bijection the mapping class group orbits of flexibility classes.

\begin{theorem}
	If $S$ is a surface of genus 2, then $\hyp(S)$ contains 9 $\MCG(S)$-orbits of flexibility classes.
\end{theorem}

\begin{proof}
	By Corollary \ref{cor:flex_class_iff_MCG-equiv}, the $\MCG(S)$-orbits of flexibility classes in $\hyp(S)$ are in bijection with signature equivalence classes of the holonomy coverings. By Theorem \ref{thm:MCG_equivelence_of_flex_covers}, there are only 12 signature equivalence classes of flexible covering by the genus 2 surface $S$. We will show that exactly 3 of these equivalence classes do not represent holonomy coverings.
	
	By Lemma \ref{lem:non-degree_1_covers}, if one of these 12 covers is not signature equivalent to a holonomy cover for some metric in $\hyp(S)$, then the signature of the covered orbifold must be $(0;2,2,2,4)$ or $(0;2,2,2,3)$. \\
	
	\noindent \textbf{Case 1: $(0;2,2,2,4)$.} Let $\orbnot$ be a $(1;2)$ orbifold and let $p_0\colon S \to O_0$ be the degree 3 flexible cover corresponding to the tuple  $$ [(012),(01), (021) ]$$ and let $q \colon O_0 \to O$ be the degree 2 branched cover induced by the involution $\tau$ shown in Figure \ref{fig:involution_torus}.  Let $\orb$ be a $(0;2,2,2,4)$ orbifold that is the quotient of $\orbnot$ by the involution  $\tau$ shown in Figure \ref{fig:involution_torus}. 
	
	\begin{figure}[ht]
		\centering
            \def\svgwidth{2.5in}
\begingroup%
  \makeatletter%
  \providecommand\color[2][]{%
    \errmessage{(Inkscape) Color is used for the text in Inkscape, but the package 'color.sty' is not loaded}%
    \renewcommand\color[2][]{}%
  }%
  \providecommand\transparent[1]{%
    \errmessage{(Inkscape) Transparency is used (non-zero) for the text in Inkscape, but the package 'transparent.sty' is not loaded}%
    \renewcommand\transparent[1]{}%
  }%
  \providecommand\rotatebox[2]{#2}%
  \newcommand*\fsize{\dimexpr\f@size pt\relax}%
  \newcommand*\lineheight[1]{\fontsize{\fsize}{#1\fsize}\selectfont}%
  \ifx\svgwidth\undefined%
    \setlength{\unitlength}{146.36205217bp}%
    \ifx\svgscale\undefined%
      \relax%
    \else%
      \setlength{\unitlength}{\unitlength * \real{\svgscale}}%
    \fi%
  \else%
    \setlength{\unitlength}{\svgwidth}%
  \fi%
  \global\let\svgwidth\undefined%
  \global\let\svgscale\undefined%
  \makeatother%
  \begin{picture}(1,0.21009545)%
    \lineheight{1}%
    \setlength\tabcolsep{0pt}%
    \put(0,0){\includegraphics[width=\unitlength,page=1]{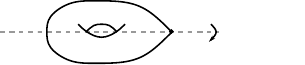}}%
    \put(0.74050351,0.10801669){\makebox(0,0)[lt]{\lineheight{1.25}\smash{\begin{tabular}[t]{l}$\tau$\end{tabular}}}}%
  \end{picture}%
\endgroup%

		\caption{The involution  of torus that gives a degree two branched cover of a $(0;2,2,2,4)$ orbifold by a $(1;2)$ orbifold.}
		\label{fig:involution_torus}
	\end{figure}
	
	Since $q$ have four branched points each with order 2, we can equip $O$ with a hyperbolic metric $\psi$ so that $(O,\psi)$ is an $(0;2,2,2,4)$ orbifold where the orbifold points  are the branch points of $q$. Pulling $\psi$ back under $q$ produces a hyperbolic metric $\psi_0$ on $O_0$ so that $\orbnot$ is a $(1;2)$ orbifold and $q$ is a locally isometric branched cover. Pulling $\psi$ back under $ q \circ p_0 $ produces a negatively curved metric $\varphi$ on $S$ so that $q \circ p_0 \colon (S,\varphi) \to \orb$ is a degree 6 flexible and $p_0 \colon (S,\varphi) \to \orbnot)$ is a degree 3 flexible cover. By Theorem \ref{thm:MCG_equivelence_of_flex_covers}, $q \circ p_0$ is signature equivalent to  the cover corresponding to the $\sym(6)$-tuple $$\bt = [ (0 2)(1 5)(3 4), (0 1)(2 4)(3 5), (0 3)(1 5)(2 4), (0 2 5 3 4 1)].$$  Since $q$ has degree larger than 1, $q\circ p_0$ cannot be signature equivalent to a holonomy covering by  Corollary \ref{cor:maximal_flexible_cover}. Hence the tuple $\bt$ does not represent a signature equivalence class of a holonomy cover. \\

	\noindent \textbf{Case 2: $(0;2,2,2,3)$.} 
	We first prove the following claim.
	
	\begin{claim}\label{claim:order_6}
		If $p \colon (S,\varphi) \to \orb$ is a degree 6 flexible cover of a $(0;2,2,2,3)$ orbifold and $p$ is not a holonomy cover, then the group generated by $\bt_p$ must have order at most 6.
	\end{claim}

	\begin{proof}
		Let $p \colon (S,\varphi) \to \orb$ be a degree 6 flexible cover of a $(0;2,2,2,3)$ orbifold. Let $r \colon (S,\varphi) \to \orbnot$ be the holonomy cover for $(S,\varphi)$.     By Corollary \ref{cor:maximal_flexible_cover}, there is a locally isometric covering map $q \colon \orbnot \to \orb$ so that $p = q \circ r$. If $q$ has degree 1, then $p$ is signature equivalent to $r$, so assume $q$ has degree at least 2. By Lemma \ref{lem:non-degree_1_covers}, this implies that $\deg(q) = 2$, $\deg(r) =3$, and $\orbnot$ has signature $(0;2,2,3,3)$. 
		
		Let $\bt_p$ be the $\sym(6)$-tuple corresponding for $p$. We  use Lemma  \ref{lem:nested_covers} to show that the group generated by $\bt_p$ must have order 6.   
		We adopt the notation from Subsection \ref{subsectnested_covers}. Let $\Sigma$ be the surface underlining $\orb$ and and $F$ be the surface underlining $\orbnot$.   Define  $\overline{\Sigma}$ and $F_0$ to be  the surfaces obtained from deleting a small neighborhood of each orbifold point of $O = \Sigma$ and $O_0 = F$ respectively. Let $\overline{F} = q^{-1}(\overline{\Sigma})$,  $\overline{S} = p^{-1}(\overline{\Sigma})$, and $S_0 = r^{-1}(F_0)$.  Note, $\Phi_{r_0}$ and $\Phi_{\overline{p}}$ are the permutation representations for $p\colon (S,\varphi)\to \orb$ and $r \colon (S,\varphi) \to \orbnot$.
		
		By Lemma \ref{lem:nested_covers}, $\Phi_{r_0}(\pi_1(F_0))$ and $\Phi_{\overline{p}}(\pi_1(\overline{F}))$ are isomorphic. Since $r$ is a flexible covering, it must be signature equivalent to one of two covers of a $(0;2,2,3,3)$ orbifold  in Theorem \ref{thm:MCG_equivelence_of_flex_covers}. Since the image of the permutation representation $\Phi_{r_0}$ is generated by the associated $\sym(3)$-tuple, we know $\Phi_{r_0}(\pi_1(F_0))$ is the order 3 cyclic subgroup of $\sym(3)$.  Thus $\Phi_{\overline{p}}(\pi_1(\overline{F}))$ is an order 3 subgroup of $\Phi_{\overline{p}}(\pi_1(\overline{\Sigma}))$. Since $\overline{q}$ is a degree 2 cover, $\Phi_{\overline{p}}(\pi_1(\overline{F}))$ has index at most 2 in $\Phi_{\overline{p}}(\pi_1(\overline{\Sigma}))$. This implies $\Phi_{\overline{p}}(\pi_1(\overline{\Sigma}))$ has order at most $6$. Hence $\bt_p$ must generate a subgroup of order at most $6$.
	\end{proof}

	Claim \ref{claim:order_6} says that the tuples 
	$$[(0 1 4 5)(2 3), (0 3)(1 5)(2 4), (0 4)(1 5)(2 3), (0 5 3)(1 2 4)]$$ and $$ [(0 3 1)(2 4 5), (0 1)(2 4)(3 5), (0 1)(2 4)(3 5), (0 1 3)(2 5 4)]$$
	must correspond to holonomy covers because they do not generate a subgroup of order at most 6. The two remaining possibilities for tuples   $$\bt_1 = [(0 3 1)(2 4 5), (0 2)(1 5)(3 4), (0 2)(1 5)(3 4), (0 1 3)(2 5 4)]$$ and 
	$$\bt_8 = [(0 3 1)(2 4 5), (0 4)(1 5)(2 3), (0 4)(1 5)(2 3), (0 1 3)(2 5 4)]$$ do generate order 6 subgroups. 
	We will now build non-holonomy covers that are signature equivalent to each of $\bt_1$ and $\bt_8$. This will show that only two of the four possible signature equivalence classes covers in this case are holonomy covers. 
	
	Fix a $(0;2,2,2,3)$ orbifold $\orb$. Let $\orbnot$ be the $(0;2,2,3,3)$ orbifold obtained by pulling back $\psi$ under the degree 2 branched cover $q$ induced by the involution in Figure \ref{fig:involution_sphere}.
	
	\begin{figure}[ht]
		\centering
            \hspace{1cm}\def\svgwidth{5in}
            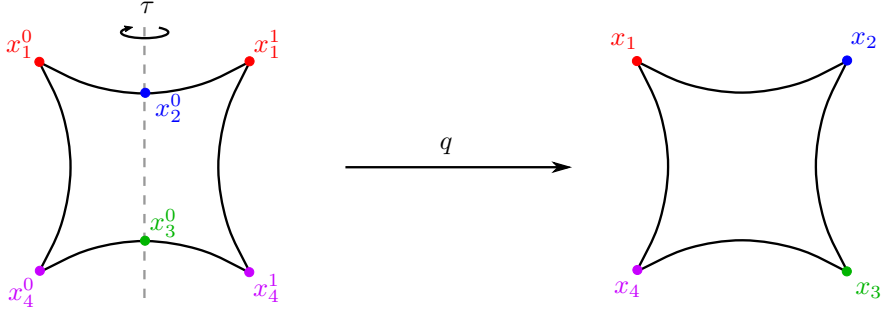
        
		\caption{The degree two covering map $q \colon \orbnot \to \orb$. The points $x_i^j$ are the preimages of $x_i$, the orbifold points of $\orb$. $x_1,x_2,x_3$ have order 2  and $x_4$ has order 3.  The order 2 orbifold points of $\orbnot$ are $x_1^0$ and $x_1^1$, while the order 3 points are $x_4^0$ and $x_4^1$. $x_2^0$ and $x_3^0$ are regular points of $\orbnot$ where the local degree of $q$ is 2.}
		\label{fig:involution_sphere}
	\end{figure}
	
	Let $r\colon S \to O_0$ be the degree three branched cover corresponding to the $\sym(3)$-tuple $$[(012),(012),(021),(021)].$$
	The composition $p=q \circ r$ is a degree 6 branched covering map from $S$ to $O$. Thus, if $\varphi$ is the pull-back of $\psi$ under $p=q \circ r$, then $r$ is a flexible covering map and $p$ is locally isometric. To see that $p$ is also a flexible cover of $\orb$ by $(S,\varphi)$, observe that the two cone points of  $(S,\varphi)$ map onto the two order 2 orbifold points of $(O_0,\psi_0)$, which then map to one of the order two orbifold points of $\orb$.
	
	We now adopt the notation of Section \ref{subsectnested_covers}: $\Sigma = O$ and $F = O_0$, then $\overline{\Sigma}$ and $F_0$ are the surfaces obtained from deleting a small neighborhood of each orbifold point of $O = \Sigma$ and $O_0 = F$ respectively. Let $\overline{F} = q^{-1}(\overline{\Sigma})$,  $\overline{S} = p^{-1}(\overline{\Sigma})$, and $S_0 = r^{-1}(F_0)$.   Note, the permutation representation $\Phi_r$ is determined by the restriction of $r$ to $S_0$ while the representation $\Phi_p$ is determined by the restriction of $p$ to $\overline{S}$. See Figure \ref{fig:lift_eta} for a diagram of the situation. 
	
		\begin{figure}
		\centering
    \def\svgwidth{6in}
    \hspace{6cm} 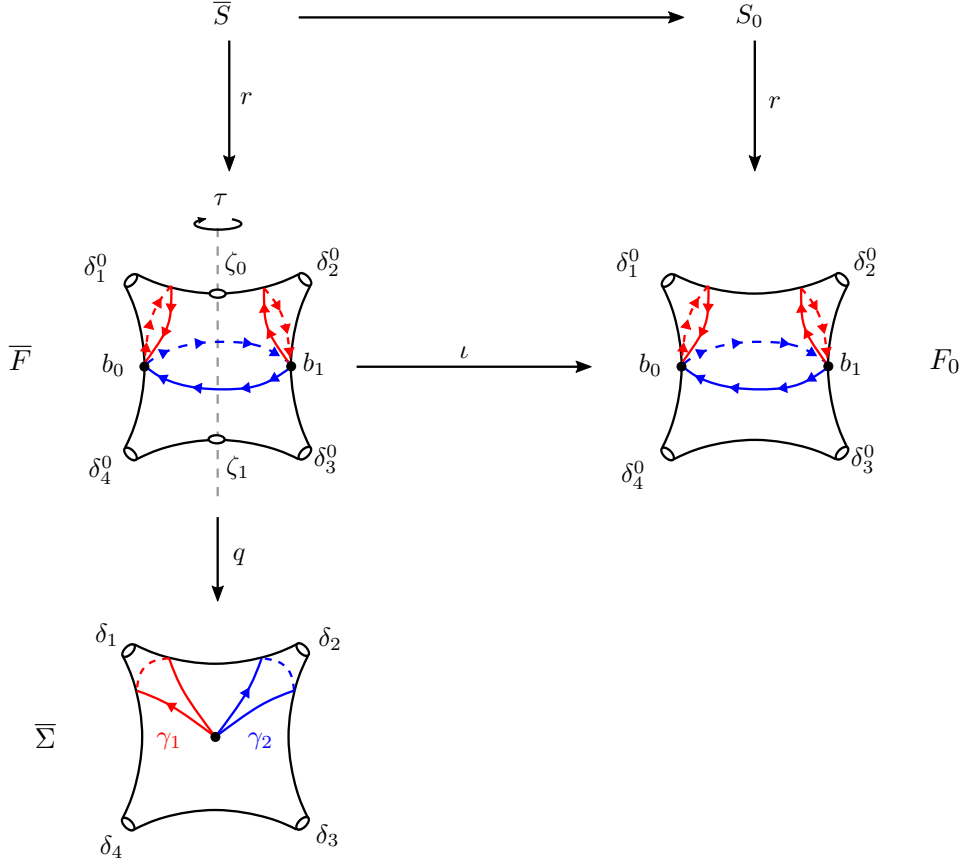
        
		\caption{The maps between $\overline{\Sigma}$, $\overline{F}$ and $F_0$. The degree 2 cover $q$ is induced  by the involution $\tau$. $\delta_1^0$ and $\delta_2^0$ map onto $\delta_1$ while $\delta_3^0$ and $\delta_4^0$ map onto $\delta_4$. $\delta_2$ and $\delta_3$ are the image of $\zeta_0$ and $\zeta_1$ respectively. $\iota$ is the inclusion of $\overline{F}$ into $F_0$. The curve $\eta = \gamma_2^{-1}\gamma_1\gamma$ lifts to a curve on $\overline{F}$ that is freely homotopic to $\delta_2^0$. }
		\label{fig:lift_eta}
	\end{figure}
	
	Let $\delta_i$ denote the boundary components of $\overline{\Sigma}$ and $\delta^0_i$ denote the boundary components of $F_0$. Now $\overline{F}$ is a subsurface of $F_0$ whose boundary contains the $\delta_i^0$, plus two additional boundary components which we call $\zeta_0$ and $\zeta_1$.   Let $a\in \overline{\Sigma}$ be the basepoint for $\pi_1(\overline{\Sigma})$ and $\{b_0,b_1\} = q^{-1}(a)$. We will let $b_0$ be the basepoint for $\pi_1(\overline{F})$ and $\pi_1(F_0)$.  As before, $\gamma_1,\gamma_2,\gamma_3,\gamma_4$ will denote the generators of $\pi_1(\overline{\Sigma})$  shown in Figure \ref{fig:4-holed_sphere}.

	Let $\bt_{p} = [s_1,s_2,s_3,s_4]$. We first argue that  the cycle structure of $s_1$  and $s_2$ are $(3,3)$ and $(2,2,2)$ respectively.  Let $x_i$ be the orbifold point of $\orb$ that are cut off by $\delta_i$.  Now $q^{-1}(x_2)$ and $q^{-1}(x_2)$ are both a single regular point of $\orbnot$. Since $r$ is flexible, this means $$p^{-1}(x_2) = r^{-1}(q^{-1}(x_2))   \text{ and } p^{-1}(x_3) = r^{-1}(q^{-1}(x_3))$$   must both be three regular point of $(S,\varphi)$  where each of the local degrees of $p$ are 2. Thus $s_2$ must have cycle structure $(2,2,2)$. Since $p$  must send the cone points of $(S,\varphi)$ to an even order orbifold point of $\orb$, we must have that $p^{-1}(x_1)$ contains both of the cone points of $(S,\varphi)$ and the local degree of $p$ at each of those cone points is 3. Thus $s_1$ has  cycle structure $(3,3)$.

	By Claim \ref{claim:order_6}, we know that  $\bt_p$ must generate an order 6 subgroup of $\sym(6)$. Consider the loop $\eta = \gamma_2^{-1} \gamma_1 \gamma_2$ in $\overline{\Sigma}$. As shown in Figure \ref{fig:lift_eta}, $\eta$ lifts to a loop $\tilde{\eta}\in \pi_1(\overline{F})$ that is freely homotopic to $\delta_2^0$. Since this loop is still freely homotopic to $\delta_2^0$ in $F_0$, $\Phi_r(\tilde{\eta}) = (012)$. Since $\gamma_1$ lifts to a loop $\tilde{\gamma_1}$ that is based at $b_0$ and freely homotopic to $\delta_1^0$ in both $\overline{F}$ and $F_0$, $\Phi_r(\tilde{\gamma_1}) = (012)$. Since $p = q \circ r$, this means the lift of  $\eta$ and $\gamma_1$ to $S$ must define the same permutation in $\sym(6)$. That is $$\Phi_p(\eta) = \Phi_{q\circ r}(\eta) = s_2^{-1}s_1s_2 = s_1 = \Phi_{q\circ r}(\gamma_1) = \Phi_p(\gamma_1).$$
	Since $s_1$ and $s_2$ have order 3 and 2 respectively, this means that $\bt_p$ must generate a cyclic subgroup of order 6. Thus, $\bt_{p}$ must by signature equivalent to    $$\bt_1 = [(0 3 1)(2 4 5), (0 2)(1 5)(3 4), (0 2)(1 5)(3 4), (0 1 3)(2 5 4)].$$

	We can repeat this argument, with $r$ being the degree 3 branched cover corresponding to $[(012),(021),(012),(021)]$ instead.  In that case, $\Phi_r(\tilde{\eta}) = (021)$ while $\Phi_r(\tilde{\gamma_1}) = (012)$. Thus, $$\Phi_p(\eta) = \Phi_{q\circ r}(\eta) = s_2^{-1}s_1s_2 \neq s_1 =\Phi_{q\circ r}(\gamma_1) = \Phi_p(\gamma_1).$$ This means that $\bt_p$ must generate a non-cyclic subgroup of order 6. So $\bt_p$ must be signature equivalent to   $$\bt_8 = [(0 3 1)(2 4 5), (0 4)(1 5)(2 3), (0 4)(1 5)(2 3), (0 1 3)(2 5 4)].$$  
\end{proof}


\section*{Code Appendix}
This appendix contains the Python functions needed to generate the conjugacy classes of $\sym(D)$-tuples $[s_1,s_2,s_3,s_4]$ where each $s_1s_2s_3 = s_4^{-1}$ and each $s_i$ has a specified cycle structure. This code is implemented using the Python intertools and sympy packages.

\begin{lstlisting}[language=Python]
from itertools import permutations
from sympy.combinatorics import Permutation
\end{lstlisting}

Care is needed when multiplying permutations with the sympy package. In this paper, the product $st$ of two permutations $s$ and $t$ is the function $s \circ t$, but in the sympy package, the product $st$ is the function $t \circ s$. That is, we would say $$(0123)(01) = (023),$$ but asking sympy to compute $(0123)(01)$ will result in $(123)$. As a result, all permutation multiplication here is written in the reserve order than the reader might expect.

First we build a Python function to collect all permutations in $\sym(D)$ with a specific cycle structure. The input `cycle\_struct' here is a list of integers for  the desired cycle structure.

\begin{lstlisting}[language=Python]
def cycle_permutations(cycle_struct, D = None):
    cycle_permutations = []
    n = sum(cycle_struct)
	
    if D == None:
        D = n
	
    #create all permutations
    all_permutations = [s  for s in permutations(set(range(D)))]
	
    #split them up according to cycle structure 
    for permutation in all_permutations:
        cycle_permutation = []
	
        last_idx = 0
        for cycle_len_idx in range(len(cycle_struct)):
            # Determine indices to cut the permutation at
            end_idx = last_idx + cycle_struct[cycle_len_idx]
            
            # add that cycle of the permuation
            cycle_permutation.append(list(permutation[last_idx: end_idx]))
            last_idx += cycle_struct[cycle_len_idx]
            
        cycle_permutations.append(Permutation(cycle_permutation))
	
    return list(set(cycle_permutations))
\end{lstlisting}

Next we have a function that checks if two $\sym(D)$-tuples $\bt_1$ and $\bt_2$ are conjugate or not. The inputs t1, t2, and Sym\_D are  lists of permutations.

\begin{lstlisting}[language=Python]
def check_sol_same(t1, t2, Sym_D):
    assert(len(t1) == len(t2))
    for sigma in Sym_D:
        satisfies = True
        for i in range(len(sol1)):
            s1 = t1[i]
            s2 = t2[i]
            
            if sigma**(-1) * s1 * sigma != s2:
                satisfies = False
                break
        if satisfies:
            return True
    return False
\end{lstlisting}

The next function takes a list of $\sym(D)$-tuples (`all\_sols') and produces a dictionary where the values are the are the conjugacy classes for the permutations  in the list and the a keys are a representative of each class.

\begin{lstlisting}[language=Python]
def get_unique_sols(all_sols, D):
    #get all group elements of Sym(D)
    Sym_D = [Permutation(s) for s in list(permutations(list(range(D))))]
	
    unique_sols = []
    equiv_classes = dict(set({})) 
	
    #store the conjugate solutions, keys are a representative for that class
    for t1 in all_sols:
        is_unique = True
        for t2 in unique_sols:
            if check_sol_same(t1, t2, Sym_D):
                is_unique = False
                if tuple(t2)in equiv_classes.keys():
                    equiv_classes[tuple(t2)].append(t1)
	
        if is_unique:
            unique_sols.append(t1)
            
            equiv_classes[tuple(t1)] = [t1]
            
    return unique_sols, equiv_classes
\end{lstlisting}

Finally we have the Python function that produces the conjugacy classes of all $\sym(D)$-tuples $[s_1,s_2,s_3,s_4]$ that satisfy $s_1s_2s_3= s_4^{-1}$ and have specified cycle structures for each $s_i$. Note, which tuple ends up as the representative for the conjugacy class is non-deterministic in this function and can change with different runs of the code.

\begin{lstlisting}[language=Python]
	
def get_sols_0_4(cycle_struct_1, cycle_struct_2,cycle_struct_3,cycle_struct_4, D):
    all_sols = []
    for s1 in cycle_permutations(cycle_struct_1, D):
        for s2 in cycle_permutations(cycle_struct_2, D):
            for s3 in cycle_permutations(cycle_struct_3, D):
            for s4 in cycle_permutations(cycle_struct_4, D):
                if s3*s2*s1 == s4**(-1):
                    all_sols.append([s1,s2,s3,s4])
    unique_sols, classes = get_unique_sols(all_sols, D)
    return unique_sols , classes
\end{lstlisting}

	\bibliography{bibliography}{}
	\bibliographystyle{alpha}
\end{document}